\documentclass[onefignum,onetabnum]{siamonline171218}

\usepackage{amsfonts}
\usepackage{graphicx}
\usepackage{dsfont}
\usepackage{subfig}
\usepackage{amssymb}
\usepackage{algorithm}
\usepackage[noend]{algpseudocode}
\ifpdf
  \DeclareGraphicsExtensions{.eps,.pdf,.png,.jpg}
\else
  \DeclareGraphicsExtensions{.eps}
\fi

\usepackage{amsmath,mathabx}
\makeatletter
\renewcommand{\pod}[1]{\allowbreak\mathchoice
  {\if@display \mkern 18mu\else \mkern 8mu\fi (#1)}
  {\if@display \mkern 18mu\else \mkern 8mu\fi (#1)}
  {\mkern4mu(#1)}
  {\mkern4mu(#1)}
}

\usepackage{enumitem}
\setlist[enumerate]{leftmargin=.5in}
\setlist[itemize]{leftmargin=.5in}

\newsiamremark{remark}{Remark}
\newsiamremark{example}{Example}
\newsiamremark{hypothesis}{Hypothesis}
\crefname{hypothesis}{Hypothesis}{Hypotheses}
\newsiamthm{claim}{Claim}
\newsiamthm{conjecture}{Conjecture}

\newcommand{\jpm}[1]{#1}

\allowdisplaybreaks

\headers{Local behavior of GEPP and GECP}{J. Peca-Medlin}

\title{Growth factors of orthogonal matrices and local behavior of Gaussian elimination with partial and complete pivoting\thanks{Submitted to the editors \today.
}
}

\author{John Peca-Medlin\thanks{Department of Mathematics, University of Arizona, Tucson, AZ 
  (\email{johnpeca@math.arizona.edu}).}
}

\usepackage{amsopn}

\makeatletter
\newcommand*{\addFileDependency}[1]{
  \typeout{(#1)}
  \@addtofilelist{#1}
  \IfFileExists{#1}{}{\typeout{No file #1.}}
}
\makeatother

\allowdisplaybreaks

\DeclareMathOperator{\diag}{diag}

\def\diag{\operatorname{diag}}
\def\V#1{{\mathbf #1}}
\def\sgn{\operatorname{sgn}}

\def\O{\operatorname{O}}
\def\P{\mathbb P}
\def\E{\mathbb E}

\def\SO{\operatorname{SO}}

\def\GL{\operatorname{GL}}

\def\Haar{\operatorname{Haar}}

\def\Uniform{\operatorname{Uniform}}

\def\Ginibre{\operatorname{Ginibre}}
\def\t{\theta}
\def\vep{\varepsilon}
\def\gepp{\operatorname{GEPP}}

\def\gecp{\operatorname{GECP}}

\ifpdf
\hypersetup{
  pdftitle={Growth factors of orthogonal matrices and local behavior of GEPP and GECP},
  pdfauthor={J. Peca-Medlin}
}
\fi

\begin{document}

\maketitle

\begin{abstract}
Gaussian elimination (GE) is the most used dense linear solver. Error analysis of GE with selected pivoting strategies on well-conditioned systems can focus on studying the behavior of growth factors. Although exponential growth is possible with GE with partial pivoting (GEPP), growth tends to stay much smaller in practice. Support for this behavior was provided \jpm{recently} by Huang and Tikhomirov's average-case analysis of GEPP, which showed GEPP growth factors \jpm{for Gaussian matrices} stay at most polynomial with very high probability. GE with complete pivoting (GECP) has also seen a lot of recent interest, with  improvements to \jpm{both lower and upper} bounds on worst-case GECP growth provided by \jpm{Bisain,} Edelman and Urschel \jpm{in 2023}. We are interested in studying how GEPP and GECP behave on the same linear systems as well as studying large growth on particular subclasses of matrices, including orthogonal matrices.  \jpm{Moreover, as a means to better address the question of why large growth is rarely encountered, we further study matrices with a large difference in growth between using GEPP and GECP, and we explore how the smaller growth strategy dominates behavior in a small neighborhood of the initial matrix.} 
\end{abstract}




\section{Introduction and motivations}

Gaussian elimination (GE) remains the most used solver of a linear system $A \V x = \V b$, where $A \in \mathbb R^{n\times n}$. GE iteratively transforms the initial matrix $A = A^{(1)}$ into an upper triangular matrix $U = A^{(n)}$ through rank one lower triangular updates to introduce zeros below each successive diagonal entry of $A^{(k)}$, the form of $A$ after $k$ steps of GE (with zeros below its first $k-1$ diagonal entries). This results in the final GE factorization $A = LU$ in $\frac23 n^3 + \mathcal O(n^2)$ total FLOPs, where $L$ is a unit lower triangular matrix. This factorization can be used to then compute solutions to the simpler triangular systems $L\V y = \V b$ and $U\V x = \V y$ using forward and backward  substitution. GE is frequently combined with a pivoting strategy, where row and column swaps are used to ensure the leading entry in the remaining untriangularized system, called the \textit{pivot}, satisfies some prescribed conditions. A pivoting strategy results instead in the factorization $PAQ = LU$ where $P$ and $Q$ are permutation matrices. This work will focus on two particular pivoting strategies, GE with partial pivoting (GEPP) and GE with complete pivoting (GECP), as well as the standard GE with no alternative pivoting strategy (GENP). (See \Cref{subsec: background} for a description of each pivoting variant.)

Understanding the behavior of GE under floating-point arithmetic has been an ongoing focus of numerical analysis for over 60 years. In exact arithmetic, GE aligns with the  $LU$ factorization of $A$, \jpm{while in floating-point arithmetic, then GE produces computed factors $\widehat L,\widehat U$ such that $A + E = \widehat L \widehat U$ for an error matrix $E$. We will adopt using $L$ and $U$ to implicitly designate \textit{computed} factors without explicitly using the notation $\widehat L, \widehat U$; this also applies for the intermediate GE forms $\widehat A^{(k)}$ as well as the associated permutation matrix factors.} Early work by Wilkinson in \cite{Wi61}, which led to the start of modern error analysis, considered studying the relative errors of {computed} solutions $\hat {\V x}$ to $A \V x = \V b$ through the bound
\begin{align}\label{eq: W ineq}
    \frac{\|\V x - \hat{\V x}\|_\infty}{\|\V x\|_\infty} \le 4n^2 \epsilon_{\operatorname{machine}} \kappa_\infty(A) \rho(A) ,
\end{align}
where $\kappa_\infty(A) = \|A\|_\infty \|A^{-1}\|_\infty$ is the $\infty$-condition number and 
\begin{align}
    \rho(A) = \jpm{\frac{\max_{i,j,k} | A_{ij}^{(k)}|}{\max_{i,j}|A_{ij}|}}
\end{align}
is the \textit{growth factor} (see \Cref{subsec: background} for other notations). In well-conditioned systems, error analysis of GE can then be reduced to analysis of the growth factor through \eqref{eq: W ineq}. 

The growth factor returns  the relative largest entry ever encountered running through the entire GE process. Using GECP, $\rho(A)$ takes the simpler form $\rho(A) = \max_k |U_{kk}|/|U_{11}|$
while in general the largest entry is not necessarily captured by  $U$. 
\jpm{Even if the largest growth is not found explicitly in this factor, it is still evident there: in \cite{Barlow_1986}, Barlow established}
\begin{equation}
	\jpm{|A_{ij}^{(k)}| \le \|L_{i,:}\|_p \|U_{:,j}\|_q}
\end{equation} 
\jpm{for $\frac1p + \frac1q = 1$ (where $L_{i,:} = \V e_i^T L$ and $U_{:,j} = U \V e_j$). Hence, using $p = \infty$ and $q = 1$ yields $\rho(A) \le \max_j \|U_{:,j}\|_1/\max_{i,j}|A_{ij}|$ when using GEPP (where $\|L_{i,:}\|_\infty = 1$ for all $i$).}

Understanding worst-case bounds on $\rho$ using different pivoting strategies has been a continued area of research since Wilkinson's initial analysis. \jpm{More recent analysis of matrix algorithm efficiency and accuracy has shifted away from worst-case analysis to more modern approaches, such as smoothed analysis (cf. \cite{Spielman_Teng_2009} and references therein), which studies  behavior under random perturbations. A full smoothed analysis using  Gaussian additive perturbations has been successfully implemented in the case of GENP growth factors by Sankar, Spielman, and Teng \cite{SST06}, but has remained out of reach for both the GEPP and GECP growth problem. The closest such result for GEPP was the recent average-case analysis work of Huang and Tikhomirov \cite{HT23}, which established high probability polynomial growth bounds using input matrices with independent and identically distributed (iid) standard normal entries, but their proof strategy cannot be upgraded to a smoothed analysis approach (i.e., they only establish bounds on  Gaussian perturbations of the zero matrix). No such (non-empirical) result for GECP has come close to smoothed analysis nor a full average-case analysis. So worst-case analysis remains relevant for these pivoting strategies.}

Our focus will be on studying the growth behavior of using both GEPP and GECP on the same linear system, and on how each strategy can inform growth behavior about the other. In particular, we are interested in local growth behavior around a system with large differences in growth behavior between both strategies. For instance, we are interested in the question of whether small perturbations on an initial system with a large discrepancy in growth behavior between both strategies concentrates toward the smaller growth for both strategies? For example, if $A$ has large GEPP growth and small GECP growth, and if $G$ is an iid Gaussian matrix, we are interested in how often $\rho^{\gepp}(A + \vep G) \approx \rho^{\gecp}(A)$ for sufficiently small $\vep$. To move toward addressing this question, we will establish bounds on how far apart differences in growth between both strategies can be when used on the same input matrix. 

We will further focus on large growth for particularly structured matrices, including orthogonal matrices, $\O(n)$ (note again we use $\mathcal O$ for standard ``big-O'' notation). This will include a refinement to the largest possible GEPP growth on orthogonal matrices, while also establishing a rich set of matrices for further study with large growth difference behavior.  Growth for structured systems has proven fruitful in recent studies. For instance, growth using GENP, GEPP and GE with rook pivoting (GERP) for \jpm{matrices formed using the Kronecker products of rotation matrices (i.e., $\bigotimes^n \SO(2)$)} is now completely understood \cite{jpm}. Although GE should not be a first choice for solving  orthogonal linear systems (viz., $Q\V x = \V b$ has the solution $\V x = Q^T \V b$ when $Q$ is orthogonal), there are situations when applying GE to orthogonal matrices makes sense. For example, Barlow needed to understand the effect of GEPP on orthogonal matrices to carry out error analysis of bidiagonal reduction \cite{B02}; this led to Barlow and Zha's analysis using GEPP on orthogonal matrices, which they showed maximized a different $L^2$-growth factor, $\rho_2$ (see \eqref{eq: L2 gf}). Additionally, while original studies of random growth factors tended to focus on ensembles with iid entries (cf. \cite{TrSc90}), many authors noted and explored the potential that orthogonal matrices can produce large growth \cite{Cryer,HiHi89,HiHi20}.  Hence, orthogonal matrices remain a rich source of study for potential large growth factors. 


\subsection{Outline of paper}
\label{subsec: overview}



Necessary notation for what follows is included in \Cref{subsec: background} while a brief history of the growth problem and preliminary results are included in \Cref{sec: prelim}. In  \Cref{sec: orth}, we  give a brief overview of the interaction of GEPP and GECP growth   for small $n$. We  also provide summary statistics for random growth samples using both Haar orthogonal and iid Gaussian matrices, which align with previous results in the literature (e.g., \cite{TrSc90}). In \Cref{subsec: large orth}, we revisit the orthogonal matrix originally studied by Barlow and Zha in \cite{BZ98} and provide  explicit intermediate GEPP forms (\Cref{prop:Q}), which provides an improved worst-case bound for $g^{\gepp}(\O(n))$ (\Cref{prop: c bound}). Proofs for results in \Cref{subsec: large orth} are postponed until \Cref{sec:thm proofs}.

In \Cref{sec: num}, we explore properties of orthogonal perturbations of linear systems. This study is analogous to the study of additive perturbations by Gaussian matrices, which was the basis for the  smoothed analysis of GENP carried out by Sankar, Speilman, and Teng as well as the recent average-case analysis of GEPP work of Huang and Tikhomirov \cite{HT23,SST06}. In these prior works, both the $L^2$-condition number, $\kappa_2$, and growth factors are studied under small additive perturbations by Gaussian matrices. Studying small multiplicative orthogonal perturbations, which preserves $\kappa_2$, enables study only of the growth factors themselves. 

In \Cref{subsec: lower bds}, we focus on addressing how much worse can GECP behave compared to GEPP in terms of growth. We employ random search paths generated through small perturbations to introduce new estimates for lower bounds for maximal GECP-GEPP growth differences for both $\GL_n(\mathbb R)$ and $\O(n)$ (see \Cref{fig:cn}). \Cref{subsec: extreme}  explores the local GEPP and GECP growth behavior for linear systems that exhibit extreme GEPP-GECP growth differences. Growth remains stable under the pivoting strategy that has minimal initial growth, while the larger growth models has local behavior progressively concentrate near the smaller initial growth. 
For our models, we then show the polynomial growth that is encountered locally limits to the initial GECP growth factor. 

\subsection{\jpm{Background and notation}} \label{subsec: background}

We will restrict focus to $n\times n$ matrices with real entries, $\mathbb R^{n\times n}$. Let $A_{ij}$ denote the element in row $i$ and column $j$ for a matrix $A \in \mathbb R^{n \times n}$. For indices $\alpha,\beta \subset [n] := \{1,2,\ldots,n\}$, let $A_{\alpha,\beta}$ denote the submatrix of $A$ built using the entries $A_{ij}$ for $i \in \alpha, j \in \beta$. Standard notation (as in MATLAB) will be used for consecutive index sequences $j:k = j,j+1,\ldots,k-1,k$, while a colon ``$:$'' is used if $\alpha = [n]$ or $\beta = [n]$ (e.g., $A_{:,j} = A \V e_j$ denotes the $j^{th}$ column of $A$, where $\V e_j$ is the standard elementary basis vector with 1 in component $j$ and 0 elsewhere). For $A \in \mathbb R^{n\times n}$, let $|A|$ denote the matrix with entries $|A|_{ij} = |A_{ij}|$, i.e., apply the absolute value entrywise to $A$. Let $\|A\|_p$ the induced matrix norm from the vector norms $\|\V x \|_p = (x_1^p + \cdots x_n^p)^{1/p}$ for $\V x \in \mathbb R^n$. Let $\V I = \sum_i \V  e_i^T \V e_i$ denote the identity matrix and $\V 0$ the zero matrix, whose dimensions are explicitly stated if not implicitly obvious. 

Let $\GL_n(\mathbb R) \subset \mathbb R^{n\times n}$ denote the group of nonsingular matrices,  $\O(n) = \{Q \in \GL_n(\mathbb R): QQ^T = \V I\}$ the subgroup of orthogonal matrices, and $\SO(n) = \{Q \in \O(n): \det Q = 1\}$ the further subgroup of special orthogonal matrices. Note this notation differs in this document from $\mathcal O(n)$, which denotes the standard big-O notation, with $f(n) = \mathcal O(g(n))$ when there exists a constant $c > 0$ such that for all sufficiently large $n$, $|f(n)| \le c|g(n)|$. Other standard complexity notation that will be used includes small-o notation ($f(n) = o(g(n))$ when$\displaystyle \lim_{n \to \infty} f(n)/g(n) = 0$) and big-$\Theta$ notation ($f(n) = \Theta(g(n))$ if $f(n) = \mathcal O(g(n))$ and $g(n) = \mathcal O (f(n))$). Let $\epsilon = \epsilon_{\operatorname{machine}}$ denote the machine epsilon, which is the smallest positive number such that $\operatorname{fl}(1+\epsilon) \ne 1$ when using floating-point arithmetic (we will use the IEEE standard notation). If using $t$-bit mantissa, then $\epsilon = 2^{-t}$. Later experiments will use double precision in MATLAB, which has 52-bit mantissa.

Standard numerical linear algebra texts (e.g., \cite{Hi02,Tr97}) can be used for a standard implementation of GENP, which iteratively transforms $A = A^{(1)}$ to $U = A^{(n)}$ using $A^{(k+1)} = (\V I + \sum_{i > k} L_{ik} \V e_i\V e_k^T)A^{(k)} = [L^{(k)}]^{-1}A$. GEPP uses only successive row pivots to ensure the pivot entry is maximal in magnitude for the leading column of the untriangularlized linear system, $A^{(k)}_{k:,k}$, while GECP uses row and column swaps to ensure the pivot maximizes the entire block $A^{(k)}_{k:,k:}$. Similar references can be used for standard implementations  of the QR factorization $A = QR$ for $A \in \GL_n(\mathbb R)$, where $Q \in \O(n)$ and $R$ is upper triangular with positive diagonal.
In particular, recall a QR factorization can be attained through the use of Householder reflectors or Givens rotations, $G(\theta,i,j)$, which is the identity matrix minimally updated so that $G(\theta,i,j)_{[i\ j],[i\ j]}$ is the (clockwise) rotation matrix with angle $\theta$. Our experiments in later sections will use Householder reflectors to efficiently sample Haar orthogonal matrices (see below), as outlined in more detail in \cite{Mezz}, along with Givens rotations to sample from an orthogonal neighborhood of a matrix (see  \cite{FrerixB19} for a similar approach).


For random variables $X,Y$, let $X \sim Y$ denote that $X$ and $Y$ are equal in distribution. Standard distributions that will be referenced include the standard Gaussian $X \sim N(\mu,\sigma^2)$ and uniform random variables $X \sim \Uniform(\mathcal S)$ for pre-compact $\mathcal S$. For $G$ a compact Hausdorff topological group, then there exists a left- and right-invariant regular probability Haar measure on $G$, $\Haar(G)$ \cite{We40}. Certain Haar and uniform measures can be sampled using standard $n \times m$ iid Gaussian models called the Ginibre ensemble, $\operatorname{Ginibre}(n,m)$. For example, for $\mathbb S^{n-1} = \{\V x \in \mathbb R^n: \|\V x\|_2 = 1\}$, if $\V x \sim \operatorname{Ginibre}(n,1)$, then $\V x/\|\V x\|_2 \sim \Uniform(\mathbb S^{n-1})$; for $B_\vep(\V 0) = \{\V x \in \mathbb R^n: \|\V x\|_2 \le \vep\}$, if $\V x \sim \operatorname{Ginibre}(n+2,1)$, then $(\vep/\|\V x\|_2)[\V I_{n} \ \V 0] \V x \sim \Uniform(B_\vep(\V 0))$; if $A \sim \operatorname{Ginibre}(n,n)$ has QR factorization $A = QR$, then $Q \sim \Haar \O(n)$. 

\section{Large growth}
\label{sec: prelim}

For nonsingular matrices $\GL_n(\mathbb R)$, we can define the maps
\begin{align}
    \rho^{\operatorname{GEPP}}: &\GL_n(\mathbb R) \to [1,a_n] \label{eq: PP interval} \\
    \rho^{\operatorname{GECP}}: &\GL_n(\mathbb R) \to [1,b_n]\label{eq: CP interval}
\end{align}
\jpm{which return the growth factors for input nonsingular matrices using the prescribed pivoting strategy.} Finding the worst-case bounds for each growth factor  is equivalent then to finding the optimal upper bounds $a_n$ and $b_n$ for each range. For $\mathcal S \subset \GL_n(\mathbb R)$, define 
\begin{align}
    g(\mathcal S) := \max_{A \in \mathcal S} \rho(A),
\end{align}
and let $g_n = g(\GL_n(\mathbb R))$, where further we use superscripts to explicitly denote a particular pivoting strategy. So optimal choices for $a_n$ and $b_n$ are then necessarily the \jpm{nondecreasing} sequences $a_n = g_n^{\operatorname{GEPP}}$ and $b_n = g_n^{\operatorname{GECP}}$.

Worst-case behavior on GEPP has been completely understood since Wilkinson's initial foray with the growth problem \cite{Wi61}: the most an entry can increase for a given GEPP step is by doubling (since $|L_{ij}| \le 1$), so the largest possible GEPP growth is attained if an entry is doubled every GEPP step, leading to the bound $\rho(A) \le 2^{n-1}$ for any $A \in \mathbb R^{n\times n}$. Moreover, Wilkinson established this bound is sharp (i.e., $g_n^{\gepp} = 2^{n-1}$) by providing a matrix in \cite[pg. 202]{Wi65}  that attains this upper bound, defined by
\begin{equation}
\label{eq:Wn}
	\jpm{W_n} = \V I_n - \sum_{i > j} \V e_i \V e_j^T + \sum_{k=1}^{n-1} \V e_n \V e_k^T,
\end{equation}
where $W_n = L_nU_n$ (in exact arithmetic) using
\begin{equation}
\label{eq: Ln}
	L_n = \V I_n - \sum_{i > j} \V e_i \V e_j^T \qquad \mbox{and} \qquad U_n = \V I_n - \V e_n\V e_n^T + \sum_{i=1}^n 2^{i-1} \V e_i \V e_n^T.
\end{equation}
$W_n$ (which needs no pivoting using GEPP) has $\max_{i,j} |W_{ij}| = 1$ and attains maximal growth (in exact arithmetic) in the corner entry  $(U_n)_{nn} = 2^{n-1}$. For example,
\begin{equation}
W_4 = \begin{bmatrix}
1 & 0 & 0& 1 \\ -1 & 1 & 0 & 1\\ -1 & -1 &1 & 1 \\ -1&-1&-1&1
\end{bmatrix} = \begin{bmatrix}
1 & 0 & 0& 0 \\ 1 & 1 & 0 & 0\\ 1 & 1 &1 & 0 \\ 1&1&1&1
\end{bmatrix}\begin{bmatrix}
1 & 0& 0& 1\\ 0&1 & 0& 2\\0&0&1&4 \\ 0&0&0&8
\end{bmatrix} = L_4U_4
\end{equation}
satisfies $\rho(W_4) = 2^3 = 8$. 

Worst-case behavior of GECP, however, remains elusive. Using Hadamard's inequality in the same work, Wilkinson established the upper bound
\begin{align} \label{eq: upper cp}
    g^{\operatorname{GECP}}_n \le (n \cdot 2\cdot  3^{1/2} \cdots n^{1/(n-1)})^{1/2} \le 2 \cdot n^{\ln n/4 + 1/2}.
\end{align}
This bound is believed to be very pessimistic. (See \cite{EdUr23} for a thorough historical and anecdotal overview of this particular conjecture as well as all relevant work on $g_n^{\operatorname{GECP}}$ up to now.) For a long time a conjecture (of apocryphal origins; cf. \cite{EdUr23}) ventured the bound $g_n^{\operatorname{GECP}}\le n$ (for real matrices).  \jpm{This conjecture survived for almost three decades}, until  Gould \cite{Go91} found a ``counterexample'' for $n = 13$ by producing a matrix using floating-point arithmetic with GECP growth of 13.0205 in 1991. Edelman \cite{EdelmanGECP} confirmed a true counterexample in exact arithmetic one year later. No substantial progress on the GECP growth problem was made for another 30 years later until last year. Edelman and Urschel \cite{EdUr23} establish a linear lower bound on $g_n^{\gecp}$, by upgrading Edelman's computational technique to upgrade Gould's floating-point ``counterexample'' into an exact arithmetic counterexample (by updating one entry by $10^{-7}$) into a theorem establishing floating-point and exact arithmetic growth factors cannot be too far from one another (cf. Theorem 4.2). \jpm{Moreover, only 6 months later, Edelman and Urschel along with now Bisain \cite{Bisain} provide the first substantial improvement to Wilkinson's upper bound for $g_n^{\gecp}$, where (using a modified Hadamard inequality) they provide a  bound of $n^{c\log n + 0.91}$ using $c \approx 0.20781$, which beats Wilkinson's original exponential $\log n$ coefficient of 0.25. The huge gulf between these lower and upper bounds leaves much on the table for future improvements.}

We are interested in how GEPP and GECP act on the same linear systems rather than just their worst-case behavior separately. That is, we are interested in establishing optimal bounds on the range of the map
\begin{align}
    \rho^{\operatorname{GEPP}} - \rho^{\operatorname{GECP}}: & \GL_n(\mathbb R) \to [-c_n,d_n] \label{eq: PP-CP interval}.
\end{align}
Studying the bounds $c_n$ and $d_n$ can then establish how much worse one pivoting strategy can behave relative to the other. 

Not a lot of research has studied the behavior between these pivoting schemes used on the same linear systems. It's clear $0 \in [-c_n,d_n]$ since if no pivoting is needed on a matrix using GECP, then that would also be the case for GEPP, so both growth factors must align; for example, this holds for $\V I$. Trivial upper bounds can be put on $c_n$ and $d_n$ using $a_n$ and $b_n$: $c_n \le b_n-1$ and $d_n \le a_n - 1 = 2^{n-1} - 1$. For $d_n$, that measures how much larger GEPP growth can be relative to GECP growth, this bound proves to be sharp.
\begin{proposition}
\begin{align}
    \max_{A \in \GL_n(\mathbb R)}(\rho^{\gepp}(A) - \rho^{\gecp}(A)) = 2^{n-1} - 1.
\end{align}
\end{proposition}
This is established by noting the matrix $\widetilde W_n = W_n(\V I_{n-1} \oplus 2)$, formed using Wilkinson's worst-case bound matrix $W_n$ (see \eqref{eq:Wn}) with its last column multiplied by 2, attains the trivial upper bound. Multiplying the last column by 2 does not change the GEPP growth factor from that for $W_n$ (see \Cref{thm: max gf}), so $\rho^{\operatorname{GEPP}}(\widetilde W_n)=\rho^{\operatorname{GEPP}}(W_n) = 2^{n-1}$, but it now reduces the GECP growth factor to the minimal value $\rho^{\operatorname{GECP}}(\widetilde W_n) = 1$: every GECP step  starts with a column swap with the last column, which is always a $\pm 2$ constant multiple of an all 1's vector while the prior columns of the $(n-k+1)\times (n-k+1)$ untriangularized block match those of $W_{n-k+1}$; this establishes no entry ever exceeds $2 = \max_{i,j}|(\widetilde W_n)_{ij}|$ at any GECP step. For example, $\widetilde W_4$ has the final GECP factorization
\begin{align}
        \begin{bmatrix}
            1&0&0&2\\-1&1&0&2\\-1&-1&1&2\\-1&-1&-1&2
        \end{bmatrix}
        \begin{bmatrix}
            0&1&0&0\\0&0&1&0\\0&0&0&1\\1&0&0&0
        \end{bmatrix}= \begin{bmatrix}
            1&0&0&0 \\1 & 1 &0&0\\ 1 & 1 & 1&0\\1&1&1&1
        \end{bmatrix}
        \begin{bmatrix}
            2 & 1&0&0\\0&-2 & 1&0\\0&0&-2&1 \\ 0&0&0&-2 
        \end{bmatrix}
    \end{align}
so that $\rho^{\operatorname{GECP}}(\widetilde W_4) = 1$.

Finding an optimal value of $c_n$ proves more interesting; after all, the trivial upper bound uses the still elusive constant $b_n = g_n^{\operatorname{GECP}}$. A na\"ive line of reasoning could lead one to assume GECP should lead to higher accuracy than GEPP and hence smaller growth, so perhaps one might expect the inequality $\rho^{\operatorname{GEPP}} \ge \rho^{\operatorname{GECP}}$. This inequality can easily be shown to hold for $n = 2$ (so the optimal interval is $[-c_2,d_2] = [0,1]$), but this fails starting at $n = 3$. This can be exhibited by
\begin{align}\label{eq: max cp-pp 3}
    B_3 = \begin{bmatrix}
        1/2 &0 &1/2\\1/2 & 1 &1 \\ 1/2 & -1 & 1
    \end{bmatrix},
\end{align}
which has $\rho^{\operatorname{GEPP}}(B_3) = 1$ and $\rho^{\operatorname{GECP}}(B_3) = 2$ (in exact arithmetic). Hence, after similarly noting $c_n$ is a nondecreasing sequence, this establishes $c_n \ge 1$ for all $n \ge 3$. Natural questions arise out of this brief exploration. \textit{How large does $c_n$ get? How far is $c_n$ from the trivial bound $b_n - 1$? How often does GEPP result in smaller growth factors than GECP? How far away can the GEPP growth factor be from the initial GECP growth factor in a small neighborhood?}

\subsection{Growth factors for small \texorpdfstring{{\boldmath$n$}}{n}}
\label{sec: orth}

While $g_n^{\gepp}$ has been known since Wilkinson's original analysis, much less is known about $g_n^{\gecp}$. The only known exact values of $g_n^{\operatorname{GECP}}$ are for $n=1$, 2, 3, $4$, which take the values exactly 1, 2, 2.25, 4. Best known lower bounds for other small $n$ are found in \cite{EdUr23}. The known constants $g_n^{\gecp}$ all carry over for $g^{\gecp}(\O(n))$ since an orthogonal matrix can be constructed to match each upper bound for $n \le 4$ (cf. \cite{Cryer,DaPe88,EdUr23}). 

\begin{figure}
    \centering
    \includegraphics[width=0.75\textwidth]{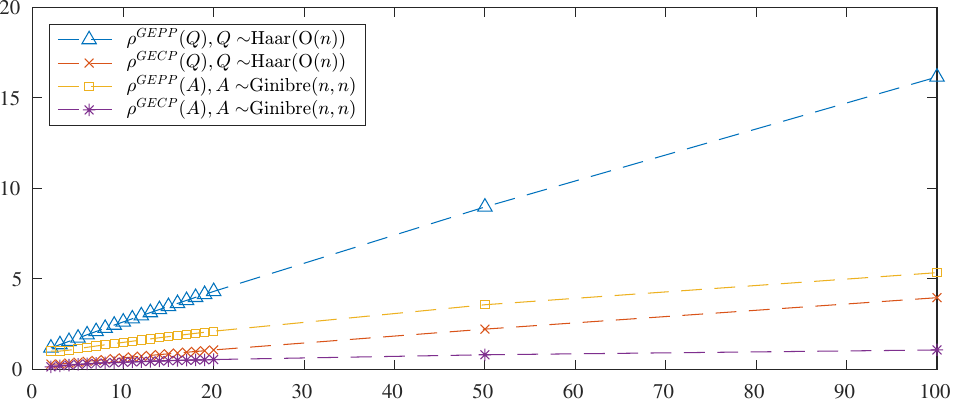}
    \caption{Plot of the average computed growth factors for $10^6$ samples of $\Haar \O(n)$ and $\Ginibre(n,n)$ for $n=2:20,50,100$.}
    \label{fig:avg gf}
\end{figure}

To better understand the overall relationship between GEPP and GECP growth factors used on the same linear systems, we will first explore this question when restricting our study to $\O(n)$. This is trivial for $n = 2$ since both growth factors necessarily align when working with $\O(n)$.
\begin{figure}[htbp]
  \centering
  \subfloat[Front view]{%
        \includegraphics[width=0.4\textwidth]{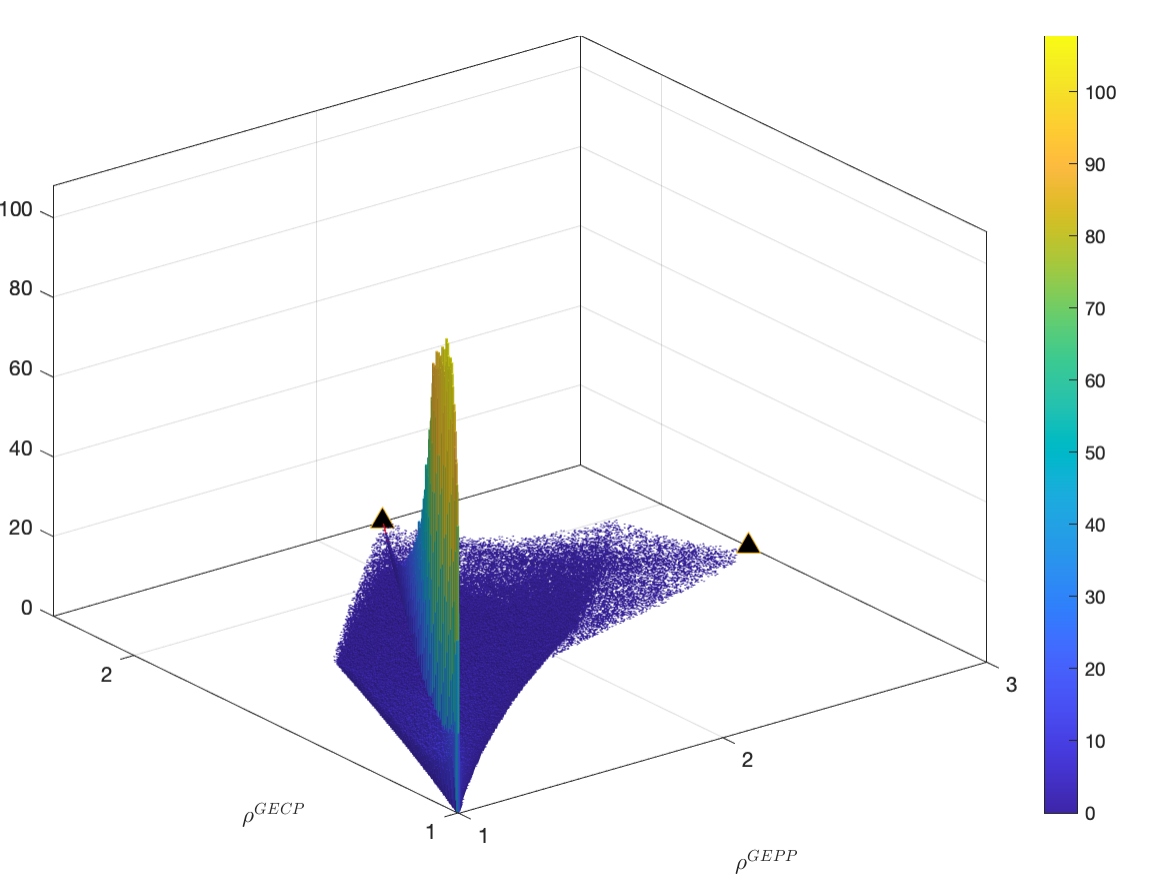}%
        }%
    \subfloat[Top view]{%
        \includegraphics[width=0.4\textwidth]{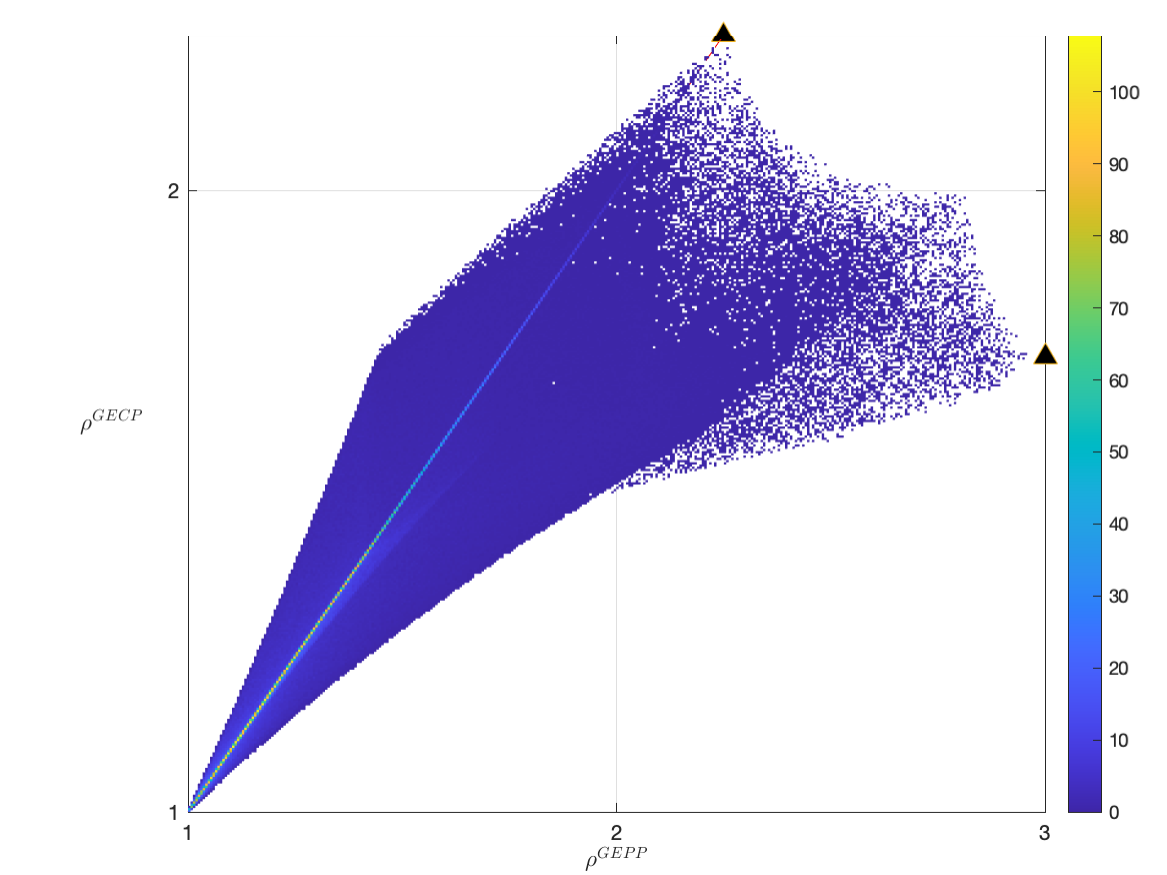}%
        }%
  \caption{Normalized histogram of $512\times 512$ grid for each pair of $\rho^{\operatorname{GEPP}}$ and $\rho^{\operatorname{GECP}}$ using $10^6$ $\Haar(\O(3))$ samples, with marks at $(2.25,2.25)$ and $(3,\sqrt3)$ corresponding to $b_3$ and $d_3$.}
  \label{fig:gf_3x3}
\end{figure}
For $n = 3$, this is no longer the case. \Cref{fig:gf_3x3} shows a top view normalized histogram (using a $512\times 512$ grid) of pairings $(\rho^{\gepp},\rho^{\gecp})$ for $10^6$ matrices sampled from $\Haar(\O(3))$. Marks are also included on the (orthogonal) $3\times 3$ matrix that attains the maximal GECP growth (of 2.25), along with the corresponding point for
\begin{align}\label{eq: max pp 3}
    Q_3 = \begin{bmatrix}
        1 & 0 & 2\\-1 & 1 & 1\\-1 & -1 & 1
    \end{bmatrix}\begin{bmatrix}
        3\\&2 \\&&6
    \end{bmatrix}^{-1/2},
\end{align}
which has $\rho^{\gepp}(Q_3) = 3$ and $\rho^{\gecp}(Q_3) = \sqrt 3$. $Q_3$ is an instance of a particular family of orthogonal matrices that will be studied more extensively in \Cref{subsec: large orth}. On \Cref{fig:gf_3x3}, $c_3$ and $d_3$ then measure how far these plotted pairs can deviate from the line $\rho^{\gepp} = \rho^{\gecp}$. \Cref{fig:avg gf} plots each average growth factor for each set of samples using $\Haar(\O(n))$ and $\Ginibre(n,n)$. Results are consistent with other studies, such as the  average-case analysis of GE carried out by Trefethen and Schreiber  \cite{TrSc90} in 1990 as well as more recent work by Higham, Higham and Pranesh \cite{HiHi20}, who establish  $\rho(Q)$ is typically not larger than $n/(4 \log n)$ for sufficiently large $n$ when $Q \sim \Haar \O(n)$ using any pivoting strategy with high probability. 

\Cref{fig:gf_3x3} exhibits high concentration for the pairs of growth factors along the line $\rho^{\gepp} = \rho^{\gecp}$. (Tables of summary statistics for $10^6$ random trials using $\Haar\O(n)$ and $\Ginibre(n,n)$ for $n=2:20,50,100$ are included in the supplementary materials in \Cref{sec: tables}.) 
For the $\Haar \O(3)$ trials, over 65\% of the samples were concentrated within $.05$ of the line $\rho^{\gepp} = \rho^{\gecp}$. The concentration along the line $\rho^{\gepp} = \rho^{\gecp}$ for $n = 3$ is more prominent if using $A \sim \Ginibre(n,n)$ (see \Cref{t:Gauss}). 

Using the sample estimates, we can empirically address the question on how often GECP growth exceeds GEPP growth for each random ensemble. Both ensembles have GEPP and GECP growth factors that have a positive fraction lie above the line $\rho^{\gecp} = \rho^{\gepp}$, i.e, a positive portion that estimates $\P(\rho^{\gecp} > \rho^{\gepp})$, with the sample proportion starting at about 13\% for $\Haar \O(3)$ and 1.4\% for $\Ginibre(3,3)$. As $n$ increases, this decreases for both ensembles, with $\Ginibre(n,n)$ in particular seeing this sample proportion decrease approximately exponentially fast, with $\P(\rho^{\gecp} > \rho^{\gepp}) \approx  e^{-n/3}$ (using a simple linear regression with the logarithmic sample proportions). 

Since $\O(n)$ is compact, then the image for the map $Q \mapsto (\rho^{\gepp}(Q),\rho^{\gecp}(Q))$ will be compact.  
For $n = 3$ at least, sampling $10^6$ $\Haar \O(3)$ matrices gives a good idea of what that image looks like (see \Cref{fig:gf_3x3}). Moreover, other attributes can be explicitly analyzed by utilizing the factorization of $\O(3)$ as the product of three Givens rotations. Using this approach, for instance, one can establish $g^{\gepp}(\O(3)) = 3$, which is achieved by $Q_3$.

\begin{figure}[htbp]
  \centering
  \subfloat[Front view]{%
        \includegraphics[width=0.4\textwidth]{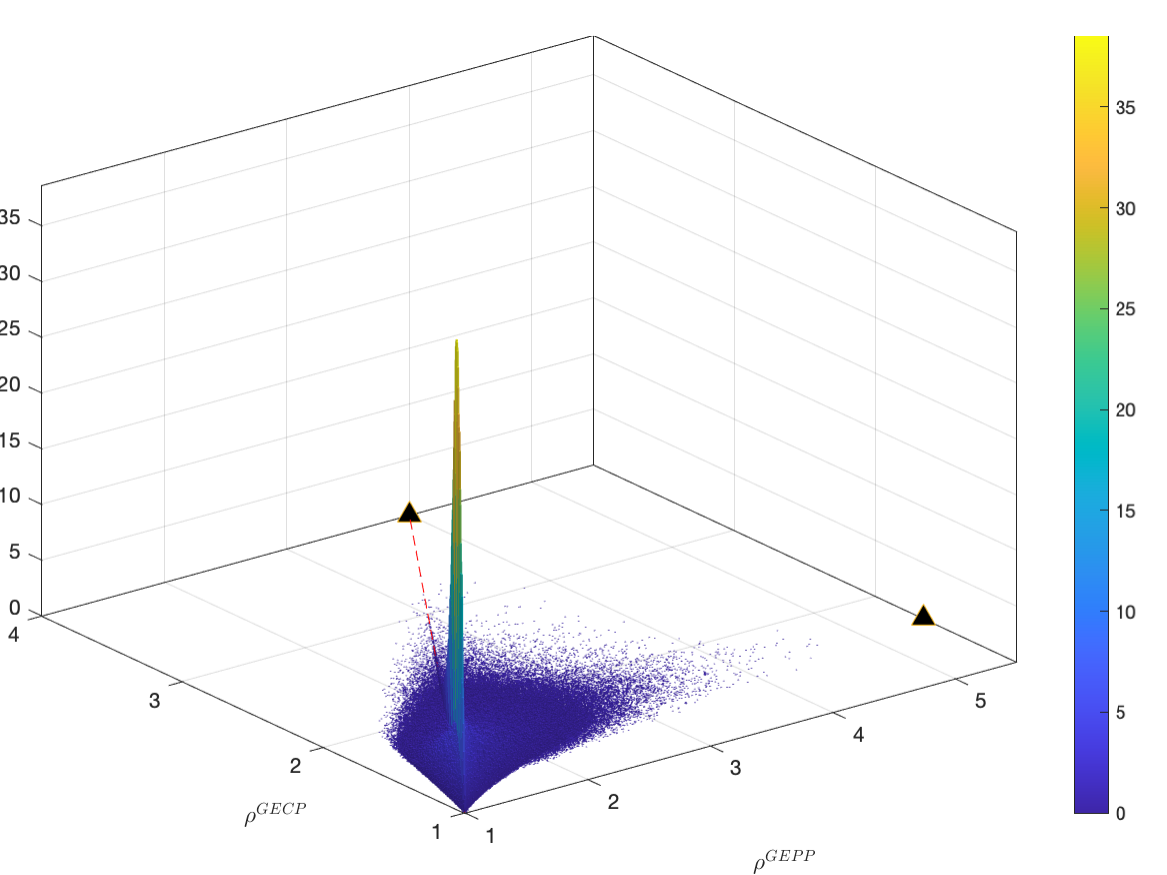}%
        }%
    \subfloat[Top view]{%
        \includegraphics[width=0.4\textwidth]{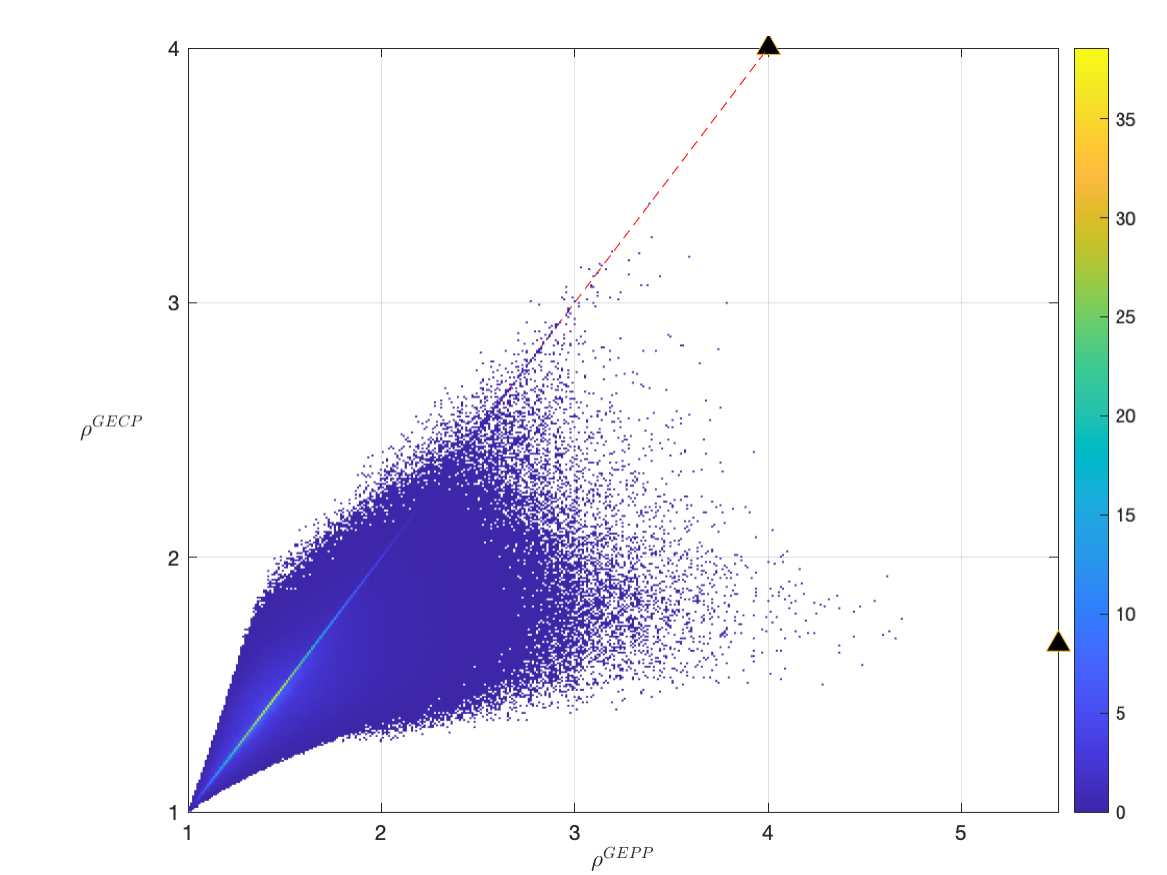}%
        }%
  \caption{Normalized histogram using $512\times 512$ grid of each pair of $\rho^{\gepp}$ and $\rho^{\gecp}$ using $10^6$ $\Haar(\O(4))$ samples, with marks at the points (4,4) (for a scaled Hadamard matrix) and $(11/2,\sqrt{11}/2)$ associated to \eqref{eq: max Q4}.}
  \label{fig:gf_4x4}
\end{figure}

This approach is untenable in general as $\dim \O(n) = {n(n-1)}/2$ provides the minimal number of Givens rotations that would be needed for such an optimization problem. For $n = 4$, uniform sampling is less effective at estimating this image, as seen in \Cref{fig:gf_4x4}. In particular, while samples were found near $(3,\sqrt 3)$ in \Cref{fig:gf_3x3}, no sample was close to the point $(11/2,\sqrt{11}/2) \approx(5.5, 1.6583)$ associated to
\begin{align}\label{eq: max Q4}
    Q_4 = \begin{bmatrix}
        1 & -1 & 0 & 4\\-1 & 5 & 0 & 2\\-1 & -3 & 1 & 1\\-1 & -3 & -1 & 1
    \end{bmatrix} \begin{bmatrix}
        4\\&44\\&&2\\&&&22
    \end{bmatrix}^{-1/2}.
\end{align}
In fact, no GEPP growth exceeded 4.7 in the $10^6$ samples, while less than $0.01\%$ had GEPP growth exceed 4. Similarly, no uniformly sampled points had GECP growth factor larger than 3.4, far from the maximal GECP growth of 4 attainable by any scaled Hadamard matrix.

From \Cref{fig:gf_3x3,fig:gf_4x4}, it appears for $\O(n)$ that the GECP growth remains much more limited in how far it can exceed GEPP growth compared to the general case. In particular, while $B_3$ from \eqref{eq: max cp-pp 3} had a difference of 1 (corresponding to the point $(1,2)$), no point for $\O(3)$ nor $\O(4)$ exhibits  coordinates differing more than  0.3	 for $n=3$ and 0.6 for $n=4$. Similarly, the extreme points associated to $a_n$, for large GEPP growth, as well as $d_n$, for how much larger GEPP growth  can be than the GECP growth, appears closest to the point for the constructed matrices we will explore in \Cref{subsec: large orth}. \Cref{sec: num} will revisit each of these constants both for the $\O(n)$ and $\GL_n(\mathbb R)$ cases.



\section{Large GEPP growth factors for \texorpdfstring{{\boldmath$\O(n)$}}{n}}
\label{subsec: large orth}

In this section, we present an explicit construction of an orthogonal matrix that attains exponential GEPP growth. We use this model to present improved estimates for the worst-case GEPP growth and max GEPP-GECP growth difference for orthogonal matrices. Proofs for results in this section are are found in \Cref{sec:thm proofs}.

\subsection{Exponential orthogonal growth model, \texorpdfstring{{$Q_n$}}{n}}

We are interested in establishing how large a GEPP growth factor can be on $\O(n)$. A trivial upper bound is $2^{n-1}$, the established maximal GEPP growth for any matrix. A perhaps natural candidate for large orthogonal growth could be $Q_n$, the orthogonal QR factor for $W_n$ from \eqref{eq:Wn}. One can compute $Q = Q_4$ in exact arithmetic (for example, applying Gram-Schmidt to the columns of $L_4$), which produces \eqref{eq: max Q4}. A straightforward check verifies $\rho^{\gepp}(Q_4) = 5.5$, which is more than 20\% larger than any of the previous $10^6$ $\Haar\O(4)$ samples encountered from \Cref{sec: orth}.

One observation, which will be formally established in \Cref{l:qr lu}, is that we would get the exact same orthogonal QR factor if we applied QR directly to $L_n$ instead of $W_n$. So the assumption that we found large orthogonal growth because we started with a matrix with a large growth  goes unfounded: $L_n$ is lower triangular and has the minimal GEPP growth factor of 1. Regardless, we will study $Q_n$, the orthogonal QR factors of $L_n$, to now establish a better bound for $g(\O(n))$.

This exact orthogonal matrix was studied extensively by Barlow and Zha in \cite{BZ98} with the goal of maximizing a different $L^2$-growth factor, defined by
\begin{align}\label{eq: L2 gf}
    \rho_2(A) := \frac{\displaystyle \max_k \|A^{(k)}\|_2}{\|A\|_2}.
\end{align}
This maximal property follows from the orthogonal invariance of the $L^2$-norm: since we can write $(QR)^{(k)} = [L^{(k)}]^{-1}QR = Q^{(k)}R$ (by \Cref{l:qr lu}), then
\begin{align}\label{eq: L2 gf ineq}
    \rho_2(QR) \le \max_k \|[L^{(k)}]^{-1}\|_2 = \rho_2(Q).
\end{align}
Since $\|A^{-1}\|_2 = 1/\sigma_{\min}(A)$, their study reduces to only establishing the smallest singular values of unit lower triangular matrices $L$ such that $|L_{ij}| \le 1$ for all $i>j$, which they show is attainable using $L_n$:
\begin{align}
    \max_{A \in \GL_n(\mathbb R)} \rho_2(A) &\le \rho_2(Q_n) = \frac1{\sigma_{\min}(L_n)} = \frac23 2^{n-1}(1 + o(1)). \label{eq:max 2 gf}
\end{align}
No explicit properties of the matrices $Q_n$ were used other than their orthogonality. Hence, using trivial norm equivalence bounds along with the $\rho_2$ maximizing properties of $Q_n$ establish 
\begin{align}
    \frac1n \cdot \frac23 2^{n-1}(1 + o(1)) &\le \rho^{\operatorname{GEPP}}(Q_n) \le n  \cdot \frac23 2^{n-1}(1 + o(1)),
\end{align}
which yields the largest GEPP growth attainable by orthogonal matrices satisfies
\begin{align}
    g^{\gepp}(\O(n)) = c 2^{n-1}(1 + o(1)) \qquad \mbox{for some} \qquad  c \in [\mbox{$\frac{2}{3n},1$}].
\end{align}
We will sharpen this bound significantly in \Cref{prop: c bound}. To do this, we will establish the explicit structure of $Q_n$ for each intermediate GEPP step.


\begin{proposition}\label{prop:Q}
Let $L_n = \V I - \sum_{i > j} \V e_i \V e_j^T$ and
\begin{align}
    \alpha_j = 1+\frac23 (4^{j-1} - 1) \quad \mbox{for $j \ge 1$}. \label{eq:alpha_j}
\end{align}
Then the QR factorization of $L_n$ has orthogonal factor given by 
    $Q_n = \hat Q \hat D^{-1/2}$, 
where
\begin{align}
    \hat Q_{:,n-1} = \V e_{n-1} - \V e_n ,\qquad
    \hat Q_{i,n} = \left\{\begin{array}{cl}
        2^{n-1 - i}, & i < n\\
        1, &i = n
    \end{array}\right. \label{eq:Q n n-1}
\end{align}
and for $n \ge 3$ and $j \le n-2$, then
\begin{align}
    \hat Q_{ij} &= \left\{\begin{array}{cl} 
    -(n-j-1)2^{j-i-1}, & {i < j}\\
    (n-j)(\alpha_j - 1) + 1, & {i = j}\\
    -\alpha_j, & i > j
    \end{array} \right. \label{eq:Qij<n-2}
\end{align}
where $\hat D$ is a diagonal matrix with nonzero entries given by
\begin{align}
    \hat D_{11} = n, \qquad \hat D_{n-1,n-1} = 2, \qquad \hat D_{nn} = 2 \alpha_{n-1}
\end{align}
and when $n \ge 4$
\begin{align}
    \hat D_{jj} &= (n-j)^2(\alpha_j - 1)^2\cdot \frac{2\alpha_j - 1}{2\alpha_j - 2} + (n-j)(\alpha_j^2+\alpha_j - 1) + 2\alpha_{j-1}, \quad 2 \le j \le n-2. \label{eq:Djj}
\end{align}

Moreover, for all $1 \le k \le n$, each intermediate GENP (and GEPP) factor takes the form
\begin{align}
    \hat Q^{(k)}_{:,n-1} = \left\{\begin{array}{cl}
        \V e_{n-1} - \V e_n, & k < n\\
        \V e_{n-1}, &k = n
    \end{array}
    \right., \quad 
    \hat Q^{(k)}_{in} = \left\{\begin{array}{cl}
        2^{n-i-1} + 2^{n-k-1}(\alpha_k - 1), & k < i < n\\
        2^{n-i-1}\alpha_i, & i \le k < n\\
        2^{n-i-1}\alpha_i, & i < k = n\\
        1 + 2^{n-k-1}(\alpha_k - 1), & k<i = n\\
        2\alpha_{n-1}, & k=i = n,
    \end{array}
    \right. \label{eq: Qk n n-1}
\end{align}
and for $n \ge 3$ then
\begin{align}\label{eq: hat Qij<n-2}
    \hat Q^{(k)}_{ij} &= \left\{\begin{array}{cl}
        -(n-j-1)[2^{j-i-1} + 2^{j-k-1}(\alpha_{k-1} - 1)], & k < i < j\\
        -(n-j-1)2^{j-i-1}\alpha_i, & i \le k < j\\
        -(n-j-1)2^{j-i-1}\alpha_i, & i < j \le k\\
        {[(n-j)(\alpha_j - 1) + 1]} - (n-j-1)2^{j-k-1}(\alpha_{k} - 1), & k < i = j\\
        \frac12(n-j+1)(\alpha_j-1) + \frac12(\alpha_j + 1), &i = j \le k\\
        -\alpha_j - (n-j-1)2^{j-k-1}(\alpha_{k} - 1), & k \le j < i\\
        0, & j < i,k.
    \end{array} \right.
\end{align}
\end{proposition}

The proof to \Cref{prop:Q}, which establishes and utilizes explicit structural properties of $Q_n$, is found in \Cref{sec:thm proofs}.

\begin{example}
    For $n = 4$, we have $Q_4$ from \eqref{eq: max Q4} is of the form $Q_4 = \hat Q \hat D^{-1/2}$ where
    \begin{align*}
        \hat Q^T \hat Q &= \begin{bmatrix}
        1 & -1 & -1 & -1\\
        -1 & 5 & -3 & -3 \\
        0 & 0 & 1 & -1\\
        4 & 2 & 1 & 2
        \end{bmatrix}
        \begin{bmatrix}
            1	& -1	& 0	& 4\\
-1&	5&	0	&2\\
-1&	-3&	1	&1\\
-1&	-3&	-1&	1
        \end{bmatrix} =
   \begin{bmatrix}
            4\\&44 \\ && 2 \\ &&&22
        \end{bmatrix} = \hat D\\
    \hat{Q}^T L_4 &= \begin{bmatrix}
        1 & -1 & -1 & -1\\
        -1 & 5 & -3 & -3 \\
        0 & 0 & 1 & -1\\
        4 & 2 & 1 & 1
        \end{bmatrix}\begin{bmatrix}
            1 & 0&0&0\\
            -1&1&0&0\\
            -1&-1&1&0\\-1&-1&-1&1
        \end{bmatrix}
        = \begin{bmatrix}
            4   &  1  &   0  &  -1\\
       &  11 &    0  &  -3\\
       &     &   2  &  -1\\
        &    &      &  1
        \end{bmatrix}\\
    \hat Q^{(2)}_{2:,2:} &= \begin{bmatrix}
        4  &   0 &    6\\
    -4   &  1  &   5\\
    -4  &  -1  &   5
    \end{bmatrix}, \quad \hat Q^{(3)}_{3:,3:}= \begin{bmatrix}
        1 &   11\\
    -1    &11
    \end{bmatrix}, \quad \hat Q^{(4)}_{44} = 22
    \end{align*} 
\end{example}

\Cref{prop:Q} can then be used to establish the GEPP growth factor of $Q_n$:
\begin{corollary}\label{cor: Q gf}
    For $Q_n$ the orthogonal QR factor of $L_n$, then
    \begin{align}
        \rho^{\operatorname{GEPP}}_n(Q_n) = \frac{2^{n-1}}{\sqrt 3}(1+o(1))
    \end{align}
\end{corollary}

\subsection{Worst-case GEPP orthogonal growth factors}

We return to exploring how large GEPP growth can be on $\O(n)$. Note first the following elementary result that highlights some of the shared properties of QR and LU factorizations. An equivalent statement is found in Proposition 1.1 in \cite{BZ98} \jpm{(the only difference is our additional focus on the signs of the diagonal entries) as well as throughout \cite{P24}}. I will include a proof for completeness.

\begin{lemma}\label{l:qr lu}
    Let $A \in \GL_n(\mathbb R)$ require no pivots using GEPP. Suppose $A$ has a QR factorization given by $A = QR$ and a GEPP factorization given by $A = LUD$ where $L$ is unit lower triangular, $U$ is upper triangular with positive diagonal, and $D$ is a diagonal sign matrix. Then $Q$ requires no pivots using GEPP while $L$ is the lower triangular factor for the GEPP factorization of $Q$ and $Q$ is the orthogonal factor for the QR factorization of $LD$.
\end{lemma}
\begin{proof}
    Each factor is nonsingular since $A$ is nonsingular while $UDR^{-1}$ and also $R (DUD)^{-1}$ are each upper triangular. Now $UDR^{-1}$ is nonsingular so no principal minors of $Q$ vanish (i.e., $Q$ has an $LU$ factorization) and so $Q$ has a unique GENP factorization given by $Q = L(UDR^{-1})$, while $Q$ requires no GEPP pivoting. Similarly, $R (DUD)^{-1}$ has positive diagonal so  $LD$ has a unique QR factorization given by $LD = Q(R(DUD)^{-1})$.
\end{proof}

\begin{remark}
    Since $W_n = L_n U_n$ requires no pivots using GEPP, then $W_n$ and $L_n$ share the same orthogonal QR factor $Q_n$ by \Cref{l:qr lu}, as previously noted. While the orthogonal factor has $L^2$-growth factor $\rho_2$ always larger than the $L^2$-growth factor for the original matrix (see \eqref{eq: L2 gf ineq}), this relationship does not hold for $\rho^{\operatorname{GEPP}}$ for $n \ge 3$, as seen by
    \begin{align}
        1 = \rho^{\operatorname{GEPP}}(L_n)  < \rho^{\operatorname{GEPP}}(Q_n) < \rho^{\operatorname{GEPP}}(W_n) = 2^{n-1}.
    \end{align}
\end{remark}

We can similarly study the orthogonal QR factors for any matrix that has maximal GEPP growth factor. Again, by \Cref{l:qr lu}, this is equivalent to studying the orthogonal QR factors for the corresponding unit lower triangular factors of a maximal GEPP growth factor matrix. Higham and Higham \cite{HiHi89} gave the explicit structure for these maximal growth models:

\begin{theorem}[\cite{HiHi89}]
\label{thm: max gf}
    If $A \in \mathbb R^{n\times n}$ has $\rho^{\operatorname{GEPP}}(A) = 2^{n-1}$, then $A = DL_n \tilde D \widecheck U \bar D$ where $D,\tilde D, \bar D$ are diagonal sign matrices, $L_n$ is as defined in \eqref{eq: Ln}, and $\widecheck U = \left[ \begin{array}{c|c}\begin{array}{c}T \\ \hline \V 0^T\end{array} & \theta \V v\end{array}\right]$ for $T$ upper triangular matrix with positive diagonal, $v_j = 2^{j-1}$ for each $j$, while $\theta = |A_{1n}| = \max_{i,j}|A_{ij}|$.
\end{theorem}

Using the trivial relationship that if $D,\tilde D$ are diagonal sign matrices then 
\begin{align}\label{eq: invar sign}
    DA^{(k)}\tilde D &= D[L^{(k)}]^{-1}A\tilde D = [(DLD)^{(k)}]^{-1}(DA\tilde D) = (DA\tilde D)^{(k)},
\end{align}
it follows  every maximal GEPP growth matrix produces an orthogonal QR factor that is sign equivalent to $Q_n$. \jpm{Let $\mathcal Q_n$ denote the orthogonal matrices of order $n$ of the form $D Q_n \tilde D$.}

\begin{corollary}\label{thm: any max Q}
    If $Q$ is the orthogonal QR factor for $A \in \mathbb R^{n\times n}$ with maximal GEPP growth, then \jpm{$Q \in \mathcal Q_n$}. Moreover, $\rho(Q) = \rho(Q_n)$ using any pivoting strategy.
\end{corollary}


\begin{remark}
    While there are uncountably many matrices with maximal growth factors, there are only finitely many such corresponding orthogonal QR factors. For $Q_n$ the orthogonal QR factor of $L_n$, then \Cref{thm: any max Q} then shows that every such orthogonal factor for a maximal GEPP growth factor matrix is of the form $Q = D Q_n\tilde D = (D Q_n D)(D \tilde D)$. Hence, \jpm{$|\mathcal Q_n| = 2^{2n-1}$} (we can assume $D_{11} = 1$). This is analogous to enumerating sign permutation equivalent Hadamard matrices. 
\end{remark}

A natural question is whether these matrices comprise the maximal possible GEPP growth among the orthogonal matrices.  Since $\mathcal Q_n \subset \O(n) \subset \GL_n(\mathbb R)$, then 
\begin{align}\label{eq:contradiction}
    \frac{2^{n-1}}{\sqrt 3}(1 + o(1)) = \rho^{\gepp}(Q_n) =  g^{\operatorname{GEPP}}(\mathcal Q_n) \le g^{\operatorname{GEPP}}(\O(n)) \le g_n^{\operatorname{GEPP}} = 2^{n-1}.
\end{align}
Moreover, when $n \ge 3$, the top inequality can be made strict: Suppose for a contradiction equality holds with $g(\O(n)) = 2^{n-1}$. By the compactness of $\O(n)$, there exists a $Q \in \O(n)$ such that $\rho^{\operatorname{GEPP}}(Q) = 2^{n-1}$. By \Cref{l:qr lu}, then $Q$ must also be the orthogonal factor for $DL_n \tilde D$ for some diagonal sign matrices $D,\tilde D$, so that then $Q \in \mathcal Q_n$ by \Cref{thm: any max Q}. Since $\max_{i,j}|Q_{ij}| \ge \|Q\V e_n\|_\infty = 2^{n-2}/\sqrt{2\alpha_{n-1}}$ \jpm{(where again $\alpha_j = 1+\frac23 (4^{j-1} - 1)$ for $j \ge 1$)}, then 
\begin{align}
    2^{n-1} = \rho^{\operatorname{GEPP}}(Q) &= \frac{\sqrt {2 \alpha_{n-1}}}{\max_{i,j}|Q_{ij}|} \le \frac{\sqrt{2\alpha_{n-1}}}{2^{n-2}/\sqrt{2\alpha_{n-1}}} 
    =2^{n-1}\cdot \frac{2 + 4^{2-n}}3< 2^{n-1},
\end{align}
 a contradiction.  This yields the following:
\begin{theorem}\label{prop: c bound}
    For $n = 1,2$, then $g^{\gepp}(\O(n)) = n$. For $n \ge 3$, there exists a constant $c \in [\frac1{\sqrt 3},1)$ such that
    \begin{align}
    g^{\operatorname{GEPP}}(\O(n)) = c2^{n-1}(1 + o(1)).
\end{align}    
\end{theorem}
 
 Finding the precise value for $c \in [\frac1{\sqrt 3},1)$ is a different ordeal. The argument used in \eqref{eq:contradiction} as well as \eqref{eq:max 2 gf} may suggest $c \in [\frac1{\sqrt 3},\frac23]$. While another reasonable guess may be that $c = \frac1{\sqrt 3}$ (we conjecture this holds, \jpm{with additional (but insufficient) support for this direction outlined in \Cref{sec: max orth prop}}). To sharpen this constant, one should maintain more of the orthogonality constraints when approaching this optimization problem, such as utilized in \cite{EAS98}. This is beyond the   scope for this paper and will be further explored in future work. 

\begin{remark}
    For comparison, empirically, GECP results in the largest pivot of $\sqrt 2$, which aligns with column $Q_n\V e_{n-1} = (\V e_{n-1} - \V e_n)/\sqrt 2$. This is easily verified for small $n$. This occurs in the final pivot for $n \le 34$, and in the antepenultimate pivot for $n \ge 35$ (at least until overflow is encountered for double precision trials). Hence, the final GECP growth factor (empirically) is
\begin{align}
    \rho^{\gecp}(Q_n) = \frac{\sqrt 2}{\max_{i,j}|(Q_n)_{ij}|} = \sqrt 2 \cdot (1 + o(1)).
\end{align}
This yields an associated lower bound on 
\begin{align}
    d_n(\O(n)) = \max_{Q \in \O(n)} \left(\rho^{\gepp}(Q) - \rho^{\gecp}(Q)\right)  \ge \left(\frac{2^{n-1}}{\sqrt 3} - \sqrt 2\right)(1 + o(1)),
\end{align}
which we further conjecture to be optimal on $\O(n)$.
\end{remark}

\subsection{GEPP Growth factors on \texorpdfstring{{\boldmath$\mathcal S_n(L)$}}{n}}

Another potential direction to tackle sharpening  bounds on the constant in \Cref{prop: c bound} is to look at how big growth can get with fixed $L$. For $L$ such that $|L_{ij}|\le 1$ for all $i > j$, let $\mathcal S_n(L)$ consist of all nonsingular matrices whose unit lower triangular GENP factor is $L$. Using the triangle inequality (for an upper bound) and then constructing a particular matrix that necessarily attains this bound, we have
\begin{align}\label{eq: max Sn}
    g^{\operatorname{GEPP}}(\mathcal S_n(L)) = 1 + \max_k \sum_{\ell = 1}^{k-1} |\gamma_{k\ell}| = \max_k \|\boldsymbol \gamma_k\|_1, \quad \gamma_{k \ell} = \sum_{\mathcal J \in \mathcal T(k,\ell)} \prod_{e(r,s) \in \mathcal J} (-L_{rs}),
\end{align}
where 
        \jpm{$\mathcal T(n,m) = \{ \{e(i_0,i_1),\ldots,e(i_{k-1},i_k)\}: n = i_0  > \cdots > i_k = m\}$ 
is the collection of decreasing paths in $\mathcal K_n$, the complete graph on $n$ vertices, connecting vertices $n > m$ (where $e(i,j)$ designates an edge connect vertices $i$ and $j$).} 
One might hope to sharpen the constant in  \Cref{prop: c bound} by establishing how close $\rho^{\gepp}(Q)$ can be to $g^{\gepp}(\mathcal S_n(L))$ for  $Q \in \O(n) \cap \mathcal S_n(L)$. This can be answered exactly for $n = 2$.
\begin{example}\label{ex: S2}
    Let $L(\alpha) = \begin{bmatrix}
        1 & 0\\\alpha & 1
    \end{bmatrix}$ for $\alpha \in [-1,1]$, which has orthogonal QR factor $Q(\alpha) = \frac1{\sqrt{1 + \alpha^2}} \begin{bmatrix}
            1 & -\alpha \\ \alpha & 1
        \end{bmatrix}$. Then $g^{\gepp}(\mathcal S_2(L(\alpha))) = 1 + |\alpha|$, which is attained by the GEPP growth factor for $\begin{bmatrix}
            1 & -\sgn(\alpha) \\ \alpha & 1
        \end{bmatrix}$, 
    while $\rho^{\gepp}(Q(\alpha)) = 1 + \alpha^2$. It follows then
    \begin{align}
        0 \le g^{\gepp}(\mathcal S_2(L(\alpha))) - \rho^{\gepp}(Q(\alpha)) = |\alpha|(1 - |\alpha|) \le \frac14,
    \end{align}
    which has equality in the lower bound iff $\alpha \in \{-1,0,1\}$ and in the upper bound iff $\alpha = \pm \frac12$.
\end{example}

Let $Q = Q(L)$ be the unique  element of $\O(n) \cap \mathcal S_n(L)$ that needs no GEPP pivots. 
Clearly $\rho^{\gepp}(Q) \le g(\mathcal S_n(L))$. Maximizing the objective function $h_{k}(\V x) = x_k + \sum_{\ell = 1}^{k-1} \gamma_{k\ell} x_\ell$ using the constraint $\V x \in \mathbb S^{n-1}$ yields the unique maximizer $\hat {\V x} = \V y/\sqrt{1 + \sum_{\ell = 1}^{k-1} \gamma_{k\ell}^2} = \V y/\|\boldsymbol \gamma_k\|_2$ where $y_i = \gamma_{ki}$ for $i < k$, $y_k = 1$ and $y_i = 0$ for $i > k$, with maximum 
    $f_k(\hat{\V x}) 
    = \|\boldsymbol \gamma_k\|_2 \le \|\boldsymbol \gamma_k\|_1 = 1 + \sum_{\ell=1}^{k-1} |\gamma_{k\ell}|$. 
This then establishes
\begin{align}\label{eq: upper bd c}
    \rho^{\gepp}(Q) \le \min\left(1, \frac{\displaystyle\max_k\|\boldsymbol \gamma_k\|_2}{\max_{i,j}|Q_{ij}| \cdot \displaystyle \max_k\|\boldsymbol \gamma_k\|_1}\right) g^{\gepp}(\mathcal S_n(L)).
\end{align}
\noindent
    This upper bound in \eqref{eq: upper bd c} is achievable (e.g., consider $L = \V I = Q$), so this does not directly sharpen the upper bound of $c$ from \Cref{prop: c bound}. Note also this inequality is off an equality by a $1 + o(1)$ factor when considering $L_n$ and $Q_n$ asymptotically, since then $\|\boldsymbol \gamma_n\|_1 = 2^{n-1}$, $\|\boldsymbol \gamma_n\|_2 = (2^{n-1}/\sqrt 3)(1 + o(1))$, and $\max_{i,j}|Q_{ij}| = 1+o(1)$.

\begin{remark}
    A future direction to potentially improve the bound on $c$ from \Cref{prop: c bound} is to establish whether  $g(\mathcal S_n(L_1)) \ge g(\mathcal S_n(L_2))$ implies $\rho(Q(L_1)) \ge \rho(Q(L_2))$ when using GEPP or whether this can fail. If this is true, this would then yield $c = 1/\sqrt 3$. This monotonicity property does hold for $n = 2$, as seen in \Cref{ex: S2}.
\end{remark}

\begin{remark}
    The idea of maximizing growth on a smaller set to attain better overall growth bounds is in lines with Theorem 3.3 in \cite{EdUr23}: if $\mathcal M \subset \mathbb R$ is bounded, then 
    \begin{align}
        g^{\gecp}(\mathcal M^{m \times m}) \ge \left(\frac{\operatorname{diam} \mathcal M}{2\operatorname{diam} |\mathcal M|} \right) g^{\gecp}_n
    \end{align}
    for all $m > 4n(3n+1)$. This attains a $n^2$ maximal growth bound factor difference for maximizing GECP growth on a restricted set, such as $\mathcal M = \{-1,0,1\}$. For GEPP, note $\{-1,0,1\}$ is sufficient to attain $g_n^{\gepp}=2^{n-1}$ as seen by $W_n$ in \eqref{eq:Wn}. 
    
    Restricting the elements of the potential orthogonal matrices, however, can either be too restrictive (e.g., only the signed permutation matrices remain if restricting $\O(n)$ to $\{-1,0,1\}$) or not restrictive enough (e.g., restricting the entries of $\O(n)$ to $[-1,1]$ returns $\O(n)$ again). Note if we restrict the elements of $\O(n)$ to $\{-\frac1{\sqrt n},\frac1{\sqrt n}\}$, then we have the scaled Hadamard matrices (when they exist, returning the empty set otherwise). Hadamard matrices have their own rich history with the growth problem (e.g., \cite{KrMi09,Wi44}).
\end{remark}

\section{GEPP and GECP growth factors in small neighborhoods}
\label{sec: num}

In this section, we will focus on the question of how much growth behavior can differ between GEPP and GECP. We will also consider the remaining worst-case growth constants when restricting focus to orthogonal matrices. From \Cref{sec: prelim}, these correspond to the constants $a_n = g_n^{\gepp}, b_n = g_n^{\gecp}$ as well as the max growth difference constants $c_n$ (measuring how much GECP growth can exceed GEPP growth) and $d_n$ (when GEPP growth exceeds GECP growth). 

\subsection{Small random perturbations}\label{subsec: perturb}
Small random perturbations are frequently used to better understand the behavior of a given cost function on a neighborhood. Additive Gaussian noise is frequently used for regularization properties, such as applied to the eigenvalue problem (e.g., \cite{BKMN21,CES22}). However, to run a similar study on the orthogonal growth problem, we can no long study behavior in additive neighborhoods, which (almost surely) move outside of $\O(n)$. Instead, we will study behavior in multiplicative orthogonal neighborhoods. These have the additional benefit of reducing the study to only the growth factors since the condition number, $\kappa_2$, is invariant under orthogonal transformations.

We use the following lemma to  stay within a desired orthogonal neighborhood of a  matrix.
\begin{lemma}\label{l: ep nb}
    Let $\boldsymbol \t = (\boldsymbol \t_1,\ldots,\boldsymbol \t_{n-1}) \in [0,2\pi)^{n(n-1)/2}$ where $\boldsymbol \t_k \in [0,2\pi)^{n-k}$. Then 
    \begin{align}
        \boldsymbol \theta \mapsto Q(\boldsymbol \t) = \prod_{i=1}^{n-1} \prod_{j = i+1}^n G(\t_{ij},i,j)
    \end{align}
    from the $2$-norm to Frobenius-normed spaces is $\sqrt{n(n-1)}$-Lipschitz continuous.
\end{lemma}
\Cref{l: ep nb} follows immediately from the orthogonal invariance of the Frobenius norm, the triangle inequality, the fact $\|U_1U_2-V_1V_2\|_F \le \|U_1 - V_1\|_F + \|U_2 + V_2\|_F$ for orthogonal $U_i,V_i$,  along with $\|\V I - G(\t,i,j)\|_F = 2\sqrt{1-\cos \t} \le \sqrt 2|\t|$. (An analogous result to \Cref{l: ep nb} is established in \cite{FrerixB19}.) This yields immediately:
\begin{corollary}\label{cor: ep nb}
    Let $A \in \GL_n(\mathbb R)$ and $\vep > 0$. If $\boldsymbol{\t}\in B_{r}(\V 0)$ for $r = \frac{\vep}{\sqrt{n(n-1)}\|A\|_F}$, then $\|A - Q(\boldsymbol{\t})A\|_F \le \vep$.
\end{corollary}
This can then be used to construct a path inside $\O(n)$ to estimate $c_n$ for $\O(n)$. 

Small Gaussian steps can similarly be achieved using standard bounds, such as:
\begin{theorem}{\cite[Theorem 4.4.5]{Ve18}}\label{thm: Ve}
    Let $A \in \mathbb R^{n\times m}$ where $A_{ij}$ are iid centered sub-Gaussian random variables. For any $t > 0$ there exists $C > 0$ such that
    \begin{align}
        \|A\| \le CK(\sqrt n + \sqrt m + t)
    \end{align}
    with very high probability, where $K = \max_{ij} \|A_{ij}\|_{\psi_2}$ for $\|X\|_{\psi_2} = \inf \{t > 0: \E e^{X^2/t^2} \le 2\}$.
\end{theorem}
    Using \Cref{thm: Ve}, one can then take small $\mathcal O(\vep)$-steps (with very high probability) inside $\GL_n(\mathbb R)$ by iteratively adding $(\vep/\sqrt n) G$ for $G \sim \Ginibre(n,n)$. We will carry this out to further estimate lower bounds for $c_n$ when using $\GL_n(\mathbb R)$.



\subsubsection{Numerical experiments}

For fixed $n$ and selected starting matrices $A \in \O(n)$, we look at the behavior of both the GEPP and GECP growth factors in small  neighborhoods of $A$. We have considered multiplicative orthogonal neighbors of the form $UA$, $AV^*$ and $UAV^*$ for $U,V$ iid $ Q(\boldsymbol \t)$ for $\boldsymbol{\t}  \sim \Uniform(B_r(\V 0))$ with $r = \frac{\vep}{\sqrt{n(n-1)}}$, but we will present only results for $UA$, which are representative of the other orthogonal  sampling methods. These will be compared to the additively perturbed matrices $A + (\vep/ \sqrt n) G$ for $G \sim \Ginibre(n,n)$. 

We will also use these sampling methods to construct random paths inside $\O(n)$ and $\GL_n(\mathbb R)$ to estimate lower bounds for the corresponding constants $c_n$, which measure how much  GECP can grow compared to GEPP, along with the other extreme range constants from \Cref{sec: prelim}. The small step sizes are much more effective with $c_n$, which appears to grow significantly slower than the other constants. Because of the small growth, we expect the lower bound to be much closer to the optimal bound for $c_n$. The other constants are better tackled using more efficient optimization software, such as JuMP (Julia for Mathematical Programming) \cite{JuMP17}, which is used in \cite{EdUr23}. Our studies are more limited by using MATLAB in double precision on a single MacBook Pro laptop, with 2 GHz Quad-Core Intel Core i5 and 16 GB 3733 MHz LPDDR4X memory. Future work will use improved methodology.

\subsection{Lower bounds for \texorpdfstring{{\boldmath$c_n$}}{n}}
\label{subsec: lower bds}

In this section, we run a random search algorithm to look for maximal difference in the GECP minus GEPP growth factors, which aligns with the constant $c_n$ from \eqref{eq: PP-CP interval} when restricting the map to $O(n)$, i.e.,
\begin{align}
    c_n(\O(n)) = \max_{Q \in \O(n)} \left(\rho^{\gecp}(Q) - \rho^{\gepp}(Q) \right).
\end{align}
For very small $n$, these estimates should be close to optimal values of $c_n$, but we would only present the remaining empirical results to just be lower bounds for larger $n$. For comparison, we also run a random search to estimate $c_n$ for all matrices, which we again strictly qualify as providing merely (sub-optimal) lower bound estimates.

The algorithm uses small random orthogonal perturbations to form a random walk inside $\O(n)$ (if the start point is inside $\O(n)$), using progressively smaller step sizes  until stopping at a matrix than attains an approximate local maximal GECP-GEPP growth difference. Pseudocode for the rudimentary random search algorithm \textsc{MaxSearch} is given in \Cref{alg: search}, where we use $f = \rho^{\gecp} - \rho^{\gepp}$. This describes a means to form a path inside $\O(n)$ that has strictly increasing GECP-GEPP growth factor differences, which stops at a point that consecutively beats $M$ of its $\vep$-neighbors. 

\begin{algorithm}
\caption{Maximal search algorithm}\label{alg: search}
\begin{algorithmic}[1]
\Procedure{MaxSearch}{$f, A \in \mathbb R^{n\times n}, \vep > 0, M$}
\State $A_0 = A$;  $j = 0$;  $k = 0$
\While {$k < M$}
\For {$\boldsymbol \t \sim \Uniform(B_{\vep}(\V 0)) \subset \mathbb R^{n(n-1)/2}$} $A_{j+1} = Q(\boldsymbol\t)A_j$
\EndFor
\If {$f(A_{j+1}) > f(A_j)$} 
\State $j \mathrel{+}= 1$;  $k = 0$
\EndIf 
\State $k \mathrel{+}= 1$
\EndWhile
\State \Return{$A_0,A_1,\ldots,A_j$}
\EndProcedure
\end{algorithmic}
\end{algorithm}

We ran several implementations of \textsc{MaxSearch} with different parameters, such as using $N=15$ starting points with initial $\vep = 10^{-1}$ and $M = 10^4$ neighbor search parameters, or $N=100$ starting points with $\vep=10^{-2}$ and $M=10^3$. Each starting point is a transformed $\Haar(\O(n))$ matrix, that uses the computed GECP permutation matrix factors to transform the starting matrix into a completely pivoted (and so partially pivoted) matrix. This results in each starting point having a growth difference of 0. The final process of refining the maximum with smaller step sizes used 9 more random searches, with parameters $\vep = 10^{-2},\ldots,10^{-10}$ and each using $M = 10^3$ neighbor comparisons. We  used an added tolerance in the comparisons of $\operatorname{tol}= 100 \epsilon_{\operatorname{machine}}$ to avoid getting stuck in a trough near $\epsilon_{\operatorname{machine}}$.

\begin{figure}
    \centering
    \includegraphics[width=0.75\textwidth]{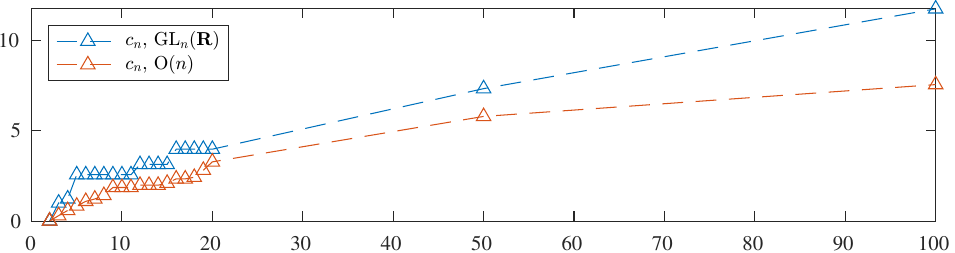}
    \caption{Plot of the computed lower bounds for $c_n$.}
    \label{fig:cn}
\end{figure}

We also compared our computations to a random search model using standard additive Gaussian perturbations. Our method differs from other implementations by initiating using $A\sim \Haar(\O(n))$ rather than $A \sim \Ginibre(n,n)$. Using $A\sim \Haar(\O(n))$ (that is then transformed to a completely pivoted matrix) leads to a starting point $(\rho^{\gepp}(A),\rho^{\gecp}(A))$ farther from the origin, which appears to be closer to the final desired extreme GECP-GEPP point (e.g., see \Cref{fig:gf_3x3,fig:gf_4x4}). Final computed lower bound estimates for each $c_n$ are displayed in \Cref{fig:cn}. (A table of the output is found in \Cref{t:c_n}.)

\subsubsection{Discussion}

This study is intended to provide a starting point for future research into how much larger can GECP growth be than GEPP growth. The particular implementation used here leaves a lot of room for improvement, such as directly using optimization software (such as JuMP), parallelization for  multiple simultaneous random search steps, as well as approximate Haar orthogonal sampling methods (such as using butterfly random orthogonal perturbations (see \cite{jpm})) for more efficient sampling. Future work will explore refinements to this initial study.

\Cref{fig:cn,t:c_n} provide some better understanding about how much worse can GECP growth be than GEPP growth. The empirical results support the lower bound estimate 
    $c_n \ge \mathcal O(n^{1/2})$ 
for both $\O(n)$ and $\GL_n(\mathbb R)$. More values of $n$ would provide a better lower bound estimate for $c_n$. In the least, we would confidently conjecture $c_n \to \infty$, i.e., GECP can always be progressively worse than GEPP. Note this lower bound estimate appears significantly far from the trivial upper bound for $c_n$ of $b_n - 1 \ge \mathcal O(n)$, where $b_n = g^{\gecp}_n$ (cf., \cite{EdUr23}). This contrasts to the result we established for the inverse question of how much worse can GEPP growth be than GECP, which we positively established matched the trivial upper bound of $2^{n-1} - 1$ for all matrices.

Also of interest would be to better establish the relationship of $c_n$ for both $\O(n)$ and $\GL_n(\mathbb R)$. Comparing the best estimates of $a_n$ for both models (i.e, $g^{\gepp}(\O(n))$ and $g_n^{\gepp}=2^{n-1}$), we established $a_n(\O(n))/a_n(\GL_n(\mathbb R) = g^{\gepp}(\O(n))/g_n^{\gepp} = c \cdot(1 + o(1))$ for some $c \in [\frac1{\sqrt 3},1)$ in \Cref{prop: c bound}. \Cref{fig:cn,t:c_n} also support the hypothesis that both constants are of the same order, with each computed ratio satisfying
\begin{align}
    c_n(\O(n))/c_n(\GL_n(\mathbb R)) \in [0.29,0.83]
\end{align}
for each $n$ in our study.

\subsection{GEPP and GECP growth factors near extreme points}
\label{subsec: extreme}

Huang and Tikhomirov \cite{HT23} establish that polynomial GEPP growth holds with high probability when using input Gaussian matrices, including small perturbations of the zero matrix. We are interested in exploring what can be said explicitly about the local growth behavior for extreme growth difference models. We will focus on analysis of a case for when GEPP growth far exceeds GECP growth here. \jpm{(An analogous case for the inverse relationship is found in \Cref{sec: CP > PP}.)}



\subsubsection{Larger GEPP than GECP growth}
\begin{figure}[htbp]
  \centering
  \subfloat[$Q(\boldsymbol \t) Q_4$]{%
        \includegraphics[width=0.4\textwidth]{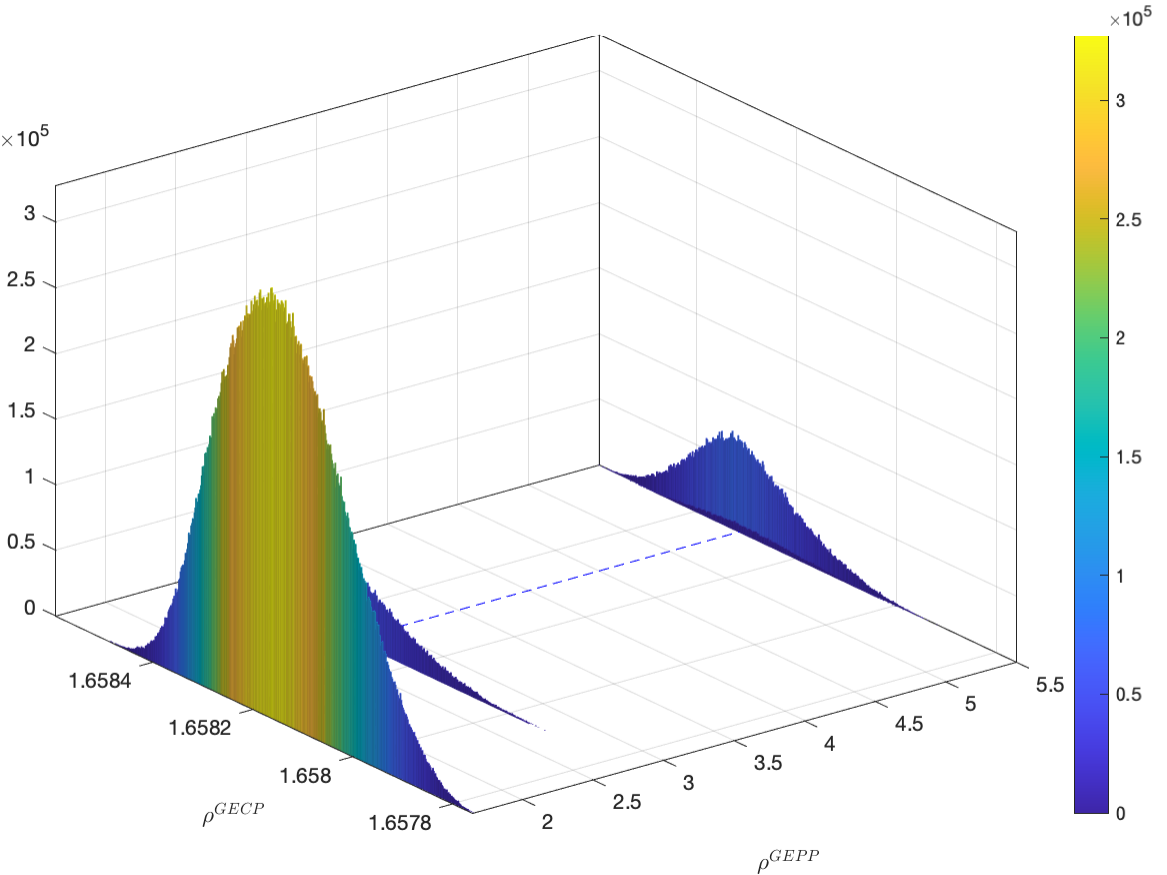}%
        }%
    \subfloat[$Q_4 + (\vep/2) G$]{%
        \includegraphics[width=0.4\textwidth]{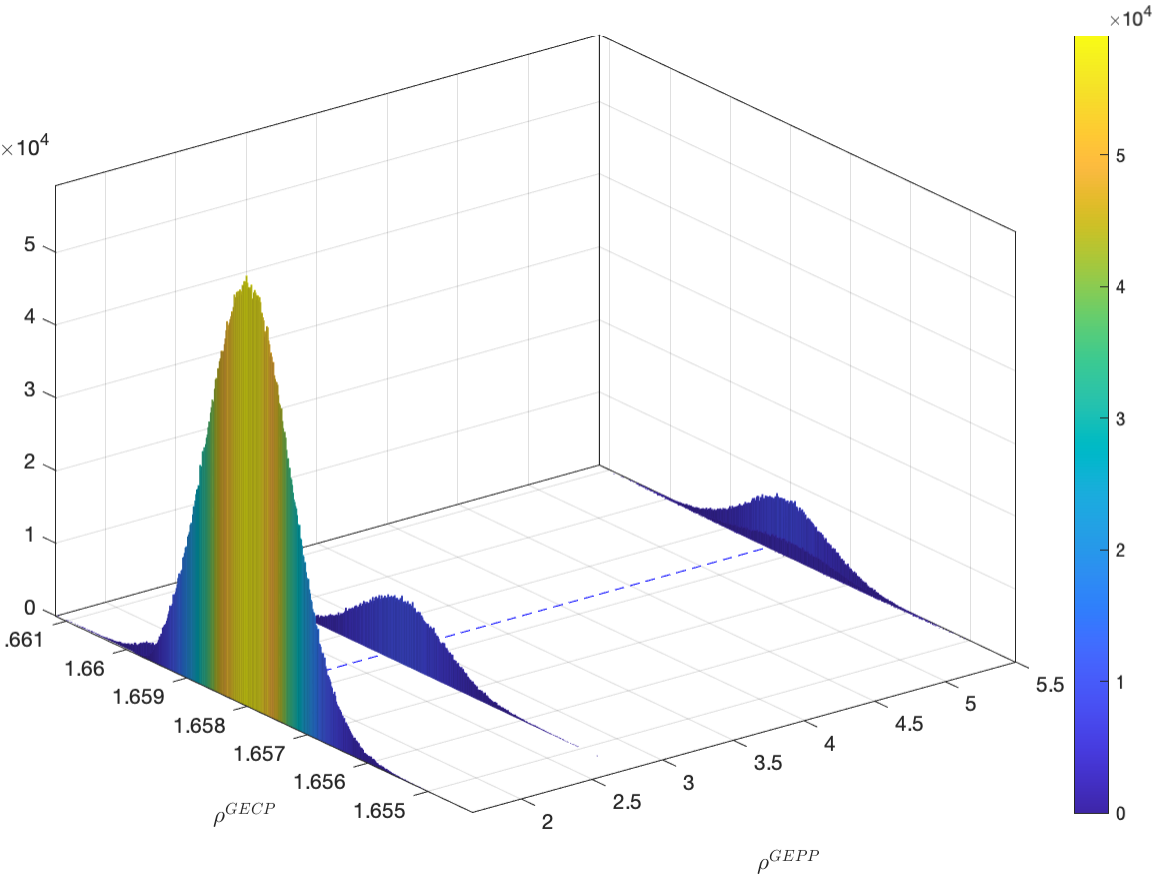}%
        }%
        \label{fig:Q4 add}
  \caption{Normalized histogram of each pair of $\rho^{\operatorname{GEPP}}$ and $\rho^{\operatorname{GECP}}$ using $10^6$ samples of (a) $UQ_4$ for $U = Q(\boldsymbol{\t})$, $\boldsymbol \t \sim \Uniform(B_r(\V 0)) \subset \mathbb R^{6}, r = \vep/\sqrt{12}$, and (b) $Q_4 + (\vep/2) G$ for $G$ for $G \sim \Ginibre(4,4)$, each using $\vep = 10^{-3}$. The line $\rho^{\gecp} = \sqrt{11}/2$ is included for comparison.}
  \label{fig:Q4 nbdh}
\end{figure}

We will next consider the local growth behavior in the case when the GEPP growth is exponential and the GECP growth moderate. \jpm{To focus the following discussion,} we will consider neighborhoods of $Q_4$, which maps to the point 
\begin{align}
    (\rho^{\gepp}(Q_4),\rho^{\gecp}(Q_4)) = (11/2,\sqrt{11},2) \approx (5.5,1.6583).
\end{align} 
$Q_n$ exhibits exponential GEPP growth of order $(2^{n-1}/\sqrt 3)(1 + o(1))$ by \Cref{cor: Q gf} and a small GECP growth factor that (empirically) concentrates on $\sqrt 2$. 
\Cref{fig:Q4 nbdh} exhibits the local GEPP and GECP growth behavior near $Q_4$, which shows a normalized histogram on a $512\times 512$ grid for $10^6$ sampled $\vep$-neighbors from both multiplicative and additive perturbations methods with $\vep = 10^{-3}$ (see \Cref{subsec: perturb}). 

The small initial GECP growth factor results in highly stable GECP growth behavior, with all $10^6$ samples remaining within $0.002$ of $\rho^{\gecp}(Q_4) = \sqrt{11}/2 \approx 1.6583$. Conversely, the local GEPP growth behavior is now highly discontinuous, with all computed $\vep$-neighbor GEPP growth factors staying within $0.002$ of either $\rho^{\gepp}(Q_4)=5.5$, $5.5/2$, or the smaller initial GECP growth factor of $1.6583$. Of these three values, most $\vep$-neighbors concentrate near the smallest limit point, with 2/3 of the additive $\vep$-neighbors having computed GEPP growth within $0.01$ of $1.6583$. This concentration near $\rho^{\gecp}(Q_n) = \sqrt 2 \cdot(1 + o(1))$ increases with $n$ \jpm{(see \Cref{t:Q_n pp-cp} for summary statistics for $\vep$-neighbors of both $Q_n$ and $W_n$)}, where $10^6$ samples for different $n$ yield the approximation 
\begin{align}
    \P\left(\rho^{\gepp}(Q_n + (\vep/\sqrt n)G) \approx \rho^{\gecp}(Q_n)\right) \approx 1 - \frac1n.
\end{align}
Similar behavior was also observed for Wilkinson's worst-case GEPP matrices, $W_n$, with again approximately $1-\frac1n$ of sampled GEPP growth factors concentrating near $\rho^{\gecp}(W_n) = 2$. Hence, for both exponential GEPP growth models $Q_n$ and $A_n$, the local GEPP growth factors using Gaussian perturbations limit to the initial GECP growth factor for each model.

For each of the extreme models studied, the smaller of the GEPP or GECP growth factors remained stable while the other exhibited highly discontinuous behavior in small neighborhoods of the initial matrix. Moreover, the local behavior progressively concentrates near the smaller initial growth factor as $n$ increases. This limiting behavior is further exhibited by each computed pair of $(\rho^{\gepp},\rho^{\gecp})$ for $\vep$-neighbors staying near the line connecting $(\rho^{\gepp}(A),\rho^{\gecp}(A))$ and $(m,m)$ for $m = \min\{\rho^{\gepp}(A),\rho^{\gecp}(A)\}$ (cf. \Cref{fig:Q4 nbdh}).



\section{Proofs of results in \Cref{subsec: large orth}}
\label{sec:thm proofs}

First we will prove \Cref{prop:Q}, which establishes the explicit intermediate GEPP forms $Q_n^{(k)}$ for $Q_n$, the orthogonal QR factor for $L_n$.

\begin{proof}[Proof of \Cref{prop:Q}]
To show $Q_n$ is orthogonal, one  needs to show $\hat Q^T\hat Q =  \hat D$, and to see the orthogonal QR factor aligns with $Q_n$, one  needs to show $\hat Q^T L_n$ or equivalently $ L_n^{-1}\hat Q$ is upper triangular with positive diagonal. This last step will also then establish each intermediate $\hat Q^{(k)}$ form as well.

A check $\hat Q_{:,j}^T \hat Q_{:,j} = \hat D_{jj}$ follows directly from \eqref{eq:Q n n-1} and \eqref{eq:Qij<n-2} along with the identity
\begin{align}\label{eq:Wn rec}
    \alpha_n + 1 = 4\alpha_{n-1}.
\end{align}
By construction $\hat Q_{:,n-1}^T \hat Q_{:,j} = 0$ for all $j \ne n-1$ is clear since $Q_{n-1:n,j}$ is a scalar multiple of $\begin{bmatrix} 1\ 1\end{bmatrix}^T$ for all $j \ne n-1$ while $Q_{:,n-1} = \V e_{n-1} - \V e_n$. Now note for $n > 2$ and $j < n-1$, \eqref{eq:Q n n-1} and \eqref{eq:Qij<n-2} can be rewritten as
\begin{align}
    \hat Q_{:,n} &= \V e_n + \sum_{k=1}^{n-1} 2^{n-1-k} \V e_k\\
    \hat Q_{:,j} &= -(n-j-1) \sum_{i=1}^{j-1} 2^{j-i-1} \V e_i + [(n-j)(\alpha_j - 1) + 1] \V e_j - \alpha_j \sum_{i = j+1}^n \V e_i.
\end{align}
It follows
\begin{align*}
    \hat Q_{:,n}^T \hat Q_{:,j} 
    &=\V e_n^T \hat Q_{:,j} + \sum_{k=1}^{j-1} 2^{n-1-k}\V e_k^T \hat Q_{:,j} + 2^{n-1-j} \V e_j^T \hat Q_{:,j} + \sum_{k=j+1}^{n-1} 2^{n-1-k}\V e_k^T \hat Q_{:,j}\\
    &= -\alpha_j -(n-j-1)2^{n-1-j}(\alpha_j-1) + 2^{n-1-j}[(n-j)(\alpha_j -1) + 1]  -\alpha_j (2^{n-1-j}-1) \\
    &=0
\end{align*}
Next, note
$
    \hat Q_{:,1} = \V e_1 - \sum_{j = 2}^n \V e_j
$ 
so that for $n \ge 4$ and $2\le j\le n-2$, then
\begin{align*}
    \hat Q_{:,1}^T \hat Q_{:,j} &= (n-j-1)\left[-2^{j-2}+\sum_{i=2}^{j-1}2^{j-i-1}\right] - [(n-j)(\alpha_j-1)+1]+\alpha_j \sum_{i=j+1}^n 1 = 0.
\end{align*}
When $n \ge 5$, an analogous and straightforward computation yields for $2 \le i < j \le n-2$
\begin{align*}
    \hat Q^T_{:,i} \hat Q_{:,j} &= (n-i-1)(n-j-1)2^{j-i-1}(\alpha_i - 1) - [(n-i)(\alpha_i - 1) + 1](n-j-1)2^{j-i-1} \\
    &\qquad + \alpha_i (n-j-1)(2^{j-i-1} - 1) - \alpha_i[(n-j)(\alpha_j - 1) + 1] + \alpha_i \alpha_j (n-j) \\
    &= 0.
\end{align*}
This establishes $\hat Q^T \hat Q = \hat D$.

Recall for $L_n^{(k)}$ the form of $L_n$ at GE step $k$, then we can write $[L^{(k)}]^{-1} \hat Q = \hat Q^{(k)}$. Recall also $L_n^{-1}$ acts on a matrix on the left by iteratively adding each row to every row below it. This last property can be reformulated as
\begin{align}
    \V e_i^T [L_n^{(k)}]^{-1} \V x &= 
        x_i + \sum_{\ell = 1}^{\min(i,k) - 1} 2^{k-1 - \ell} x_\ell,
\end{align}
for any $\V x \in \mathbb R^{n}$. Using  \eqref{eq:Q n n-1} and the fact 
\begin{align}
    \sum_{\ell = 1}^{m-1} 2^{k-1-\ell}2^{n-\ell - 1} = 2^{n+k-2m-1}(\alpha_m - 1),
\end{align}
yields then directly  \eqref{eq: Qk n n-1}, using also \eqref{eq:Wn rec} for the $k = i = n$ case along with $1 + \frac12(\alpha_n - 1) = 2\alpha_{n-1}$. \eqref{eq: Qk n n-1} follows even more directly from \eqref{eq:Q n n-1}. 

Since 
\begin{align}
    \hat Q_{1:j-1,j} = -(n-j-1) \begin{bmatrix}
        2^{j-i-1}\ \cdots \ 2 \ 1
    \end{bmatrix}^{T} = -(n-j-1) \hat Q_{n-j+1:n-1,n}
\end{align}
we can reuse  \eqref{eq: Qk n n-1}, the intermediate forms for the last column of $Q$,  to then yield \eqref{eq: hat Qij<n-2}; note also $\hat Q_{ij}^{(k)} = 0$ for $k > j$ and $i > j$ since for $k = j$, then
\begin{align}
    \hat Q_{jj}^{(j)} 
    &= \frac12(n-j+1)(\alpha_j-1) + \frac12(\alpha_j + 1) \\
    \hat Q_{ij}^{(j)} &= -\left[\alpha_j + (n-j-1)2^{-1}(\alpha_{j} - 1) \right] = -\hat Q_{jj}^{(j)}
\end{align}
so the following $(j+1)$-GE step  eliminates all entries below the $jj$ entry. Since  $0 = \hat Q_{ij}^{(j+1)} =\hat Q_{ij}^{(j+2)} = \cdots =  \hat Q_{ij}^{(n)}$, this establishes $L_n^{-1} \hat Q = \hat Q^{(n)}$ is upper triangular. Moreover, since \begin{align}
    \hat Q_{jj}^{(n)} = 1 + \frac12(n-j+1)(\alpha_j - 1) \mbox{ for } j < n-1, \quad 
    \hat Q_{n-1,n-1}^{(n)} = 1,\quad
    \hat Q_{nn}^{(n)} = 2\alpha_{n-1},
\end{align}
then $\hat Q^{(n)}$ has positive diagonal.
\end{proof}

Next, we will use \Cref{prop:Q} to establish the asymptotic GEPP growth factor for $Q_n$.

\begin{proof}[Proof of \Cref{cor: Q gf}]
    Using \eqref{eq: hat Qij<n-2} for $2\le j \le n-2$, we have each GENP (and GEPP) intermediate form $\hat Q^{(k)}$ satisfies
    \begin{align*}
        \begin{array}{c}\displaystyle \max_{k}\|\hat Q_{1:j-1,j}^{(k)}\|_\infty = |\hat Q_{j-1,j}^{(n)}| = (n-j-1)\alpha_{j-1}, \  
         \max_{k}\|\hat Q_{jj}^{(k)}\|_\infty = |\hat Q_{jj}| = (n-j)(\alpha_j - 1) + 1\\
        \displaystyle\max_k \|\hat Q_{j+1:n,j}^{(k)}\|_\infty = |\hat Q^{(j)}_{j+1,j}| = \frac12(n-j+1)(\alpha_j-1) + 2 \alpha_{j-1}.
        \end{array}
    \end{align*}
    Since also 
        $\max_k \|\hat Q^{(k)}_{:,1}\|_\infty = |\hat Q_{11}| =  1 = |\hat Q_{n-1,n-1}| = \max_k \|\hat Q^{(k)}_{:,n-1}\|_\infty$ 
    then 
        $\max_k \|\hat Q^{(k)}_{:,j}\|_{\max} = |\hat Q_{jj}|$ for all $j \le n-1$. 
    Together, this shows
    \begin{align}
        \max_k \|Q^{(k)}_{:,1:n-1}\|_{\max} = \max_k \{|\hat Q_{jj}|D_{jj}^{-1/2}: \ j=1,\ldots,n-1\} \le \|Q\|_{\max} \le \|Q\|_2 = 1,
    \end{align}
    while by \eqref{eq: Qk n n-1} we have 
        $\max_k \|Q^{(k)}_{:,n}\|_\infty = |\hat Q^{(n)}_{nn}|D^{-1/2}_{nn} = 2\alpha_{n-1}D_{nn}^{-1/2} = \sqrt{2\alpha_{n-1}} > 1$. 
    This yields 
        $\max_{i,j,k} |Q^{(k)}_{ij}| = \sqrt{2\alpha_{n-1}}$. 
    Moreover, by \Cref{prop:Q}, we have
    \begin{align*}
        |Q_{jj}| &= |\hat Q_{jj}| D_{jj}^{-1/2} 
        = \frac{(n-j)(\alpha_j -1) + 1}{\left[(n-j)^2(\alpha_j - 1)^2\cdot \frac{2\alpha_j - 1}{2\alpha_j - 2} + (n-j)(\alpha_j^2+\alpha_j - 1) + 2\alpha_{j-1}\right]^{1/2}}&\\
        &= \sqrt{\frac{2\alpha_j - 2}{2 \alpha_j -1}}\frac{1 +  \Theta\left(\frac1{(n-j)\alpha_j}\right)}{\left[1 + \Theta(\frac1{n-j})+\Theta(\frac1{(n-j)\alpha_j}) \right]^{1/2}}
    \end{align*}
    for each $j \le n-2$. Taking $j = j(n)$ so that $n - j = \omega_n(1)$, 
    then 
        $1 + o(1) = |Q_{j(n),j(n)}| \le \max_{i,j}|Q_{ij}| \le 1$  {and so also } $1/{\max_{i,j}|Q_{ij}|} = 1+o(1)$. 
%
    Note 
    \begin{align}\label{eq: max Qk asym}
        \sqrt{2\alpha_{n-1}} = \left[2 \left(1 + \frac23(4^{n-2} - 1)\right)\right]^{1/2} = \frac{2^{n-1}}{\sqrt 3}\sqrt{1 + 2^{3-2n}} = \frac{2^{n-1}}{\sqrt 3}(1 + o(1)).
    \end{align}
    Combining these along with the fact $|L_{ij}|\le 1$ for all $i > j$ yields
    \begin{align}
        \rho^{\operatorname{GEPP}}(Q) = \frac{2^{n-1}}{\sqrt 3}(1 + o(1)).
    \end{align}
    
    
\end{proof}

\appendix

\section{Supporting tables}
\label{sec: tables}

Tables containing the summary data for the trials or experiments used in \Cref{sec: prelim}, \Cref{subsec: lower bds}, and \Cref{subsec: extreme} are included here. 

\Cref{t:O(n),t:Gauss} consists of the summary statistics for the $10^6$ trials using both $\Haar \O(n)$ and $\Ginibre(n,n)$ matrices for $n = 2:20,50,100$, which  includes the sample medians, the sample means ($\bar x$), and the sample standard deviations ($s$). The last three columns include sample proportions that, respectively, satisfied the inequalities $\rho^{\gepp} + \operatorname{tol}< \rho^{\gecp}$, $|\rho^{\gepp} - \rho^{\gecp}| \le \operatorname{tol}$, and $\rho^{\gepp} - \operatorname{tol} >\rho^{\gecp}$ where $\operatorname{tol} = 0.05$. These respective sample proportions are used to estimate $\P(\rho^{\gepp} < \rho^{\gecp})$, $\P(\rho^{\gepp} = \rho^{\gecp})$, and $\P(\rho^{\gepp} > \rho^{\gecp})$.

\Cref{t:c_n} consists of the resulting best computed lower bounds of $c_n$ for both $\O(n)$ and $\GL_n(\mathbb R)$. Descriptions of the computation methods are outlined in \Cref{subsec: lower bds}.

\Cref{t:Q_n pp-cp} consists of summary statistics on $10^6$ sampled $\vep$-neighborhoods of $Q_n$ and $W_n$ for $n = 2:20,50,1000$, using $\vep = 10^{-3}$. The final column for each subsection includes the sample proportion to estimate $\P(X_n = 0)$ and $\P(Y_n = 0)$, where $X_n = \rho^{\gecp}(Q_n) - \rho^{\gepp}(Q_n + \overline G)$ and $Y_n = \rho^{\gecp}(W_n) - \rho^{\gepp}(W_n + \overline G)$ for $\overline G = (\vep/\sqrt n) G$ with $G \sim \Ginibre(n,n)$. An additional tolerance of $\operatorname{tol} = 0.01$ was use, so that $\P(|Z| < \operatorname{tol}) \approx \P(X = 0)$ for $Z = X_n$ or $Z = Y_n$.

\begin{table}[!htbp]
\centering
{
\begin{tabular}{r|ccc|ccc|ccc}
     &\multicolumn{3}{c}{$\rho^{\gepp}$} &\multicolumn{3}{|c}{$\rho^{\gecp}$} &\multicolumn{3}{|c}{$\approx \P(\rho^{\gepp}  (<,=,>)  \rho^{\gecp})$}\\
     $n$ & Median & $\bar x$ & $s$& Median &  $\bar x$ & $s$ & $<$ & $=$ & $>$\\ \hline  
    2	&	1.1714	&	1.2730	&	0.2765	&	1.1714	&	1.2730	&	0.2765	&	-	&	1	&	-	\\
3	&	1.3571	&	1.4110	&	0.2847	&	1.3557	&	1.3795	&	0.2223	&	0.1297	&	0.6536	&	0.2167	\\
4	&	1.5372	&	1.6061	&	0.3381	&	1.4949	&	1.5200	&	0.2256	&	0.1696	&	0.4355	&	0.3949	\\
5	&	1.7191	&	1.7917	&	0.3857	&	1.6040	&	1.6308	&	0.2244	&	0.1545	&	0.3058	&	0.5397	\\
6	&	1.9017	&	1.9807	&	0.4357	&	1.7085	&	1.7362	&	0.2330	&	0.1290	&	0.2216	&	0.6494	\\
7	&	2.0828	&	2.1678	&	0.4832	&	1.8098	&	1.8371	&	0.2429	&	0.1048	&	0.1652	&	0.7300	\\
8	&	2.2606	&	2.3521	&	0.5293	&	1.9085	&	1.9355	&	0.2537	&	0.0856	&	0.1245	&	0.7899	\\
9	&	2.4372	&	2.5382	&	0.5780	&	2.0072	&	2.0336	&	0.2650	&	0.0711	&	0.0955	&	0.8335	\\
10	&	2.6118	&	2.7212	&	0.6238	&	2.1032	&	2.1295	&	0.2762	&	0.0593	&	0.0743	&	0.8664	\\
11	&	2.7858	&	2.9042	&	0.6707	&	2.1984	&	2.2242	&	0.2867	&	0.0497	&	0.0591	&	0.8912	\\
12	&	2.9570	&	3.0837	&	0.7150	&	2.2918	&	2.3173	&	0.2972	&	0.0428	&	0.0475	&	0.9097	\\
13	&	3.1283	&	3.2633	&	0.7598	&	2.3839	&	2.4097	&	0.3077	&	0.0370	&	0.0382	&	0.9248	\\
14	&	3.2987	&	3.4430	&	0.8048	&	2.4754	&	2.5013	&	0.3185	&	0.0316	&	0.0317	&	0.9368	\\
15	&	3.4679	&	3.6196	&	0.8473	&	2.5652	&	2.5909	&	0.3280	&	0.0278	&	0.0260	&	0.9461	\\
16	&	3.6355	&	3.7957	&	0.8919	&	2.6546	&	2.6802	&	0.3376	&	0.0240	&	0.0218	&	0.9541	\\
17	&	3.8013	&	3.9707	&	0.9333	&	2.7430	&	2.7694	&	0.3480	&	0.0214	&	0.0184	&	0.9603	\\
18	&	3.9709	&	4.1474	&	0.9774	&	2.8303	&	2.8569	&	0.3572	&	0.0186	&	0.0155	&	0.9660	\\
19	&	4.1331	&	4.3168	&	1.0156	&	2.9156	&	2.9424	&	0.3672	&	0.0166	&	0.0135	&	0.9698	\\
20	&	4.2981	&	4.4928	&	1.0625	&	3.0020	&	3.0285	&	0.3764	&	0.0149	&	0.0114	&	0.9737	\\
50	&	8.9687	&	9.3818	&	2.2129	&	5.3573	&	5.4007	&	0.6380	&	0.0016	&	0.0007	&	0.9978	\\
100	&	16.156	&	16.9045	&	3.9513	&	8.8896	&	8.9626	&	1.0245	&	0.0003	&	0.0001	&	0.9996	
\end{tabular}
}
\caption{Summary statistics for computed $\rho^{\gepp}(A)$ and $\rho^{\gecp}(A)$ using $10^6$ samples  $A \sim \Haar\O(n)$}
\label{t:O(n)}
\end{table}

\begin{table}[!htbp]
\centering
{
\begin{tabular}{r|ccc|ccc|ccccc}
     &\multicolumn{3}{c}{$\rho^{\gepp}$} &\multicolumn{3}{|c}{$\rho^{\gecp}$} &\multicolumn{3}{|c}{$\approx \P(\rho^{\gepp}  (<,=,>)  \rho^{\gecp})$}\\
     $n$ & Median & $\bar x$ & $s$& Median &  $\bar x$ & $s$ & $<$ & $=$ & $>$\\ \hline 
 2	&	1.0000	&	1.0438	&	0.1247	&	1.0000	&	1.0112	&	0.0598	&	-	&	0.871459	&	0.128541	\\
3	&	1.0000	&	1.0974	&	0.1811	&	1.0000	&	1.0270	&	0.0883	&	0.0139	&	0.7060	&	0.2802	\\
4	&	1.0480	&	1.1549	&	0.2235	&	1.0000	&	1.0415	&	0.1038	&	0.0255	&	0.5494	&	0.4250	\\
5	&	1.1315	&	1.2167	&	0.2602	&	1.0000	&	1.0558	&	0.1153	&	0.0317	&	0.4144	&	0.5539	\\
6	&	1.2076	&	1.2804	&	0.2905	&	1.0000	&	1.0694	&	0.1236	&	0.0326	&	0.3067	&	0.6607	\\
7	&	1.2789	&	1.3466	&	0.3184	&	1.0000	&	1.0829	&	0.1302	&	0.0303	&	0.2214	&	0.7483	\\
8	&	1.3463	&	1.4131	&	0.3422	&	1.0220	&	1.0959	&	0.1354	&	0.0260	&	0.1572	&	0.8168	\\
9	&	1.4136	&	1.4811	&	0.3642	&	1.0491	&	1.1094	&	0.1401	&	0.0209	&	0.1104	&	0.8687	\\
10	&	1.4803	&	1.5489	&	0.3839	&	1.0735	&	1.1231	&	0.1445	&	0.0165	&	0.0766	&	0.9070	\\
11	&	1.5447	&	1.6160	&	0.4026	&	1.0962	&	1.1370	&	0.1482	&	0.0124	&	0.0532	&	0.9345	\\
12	&	1.6090	&	1.6827	&	0.4197	&	1.1175	&	1.1511	&	0.1515	&	0.0095	&	0.0359	&	0.9546	\\
13	&	1.6723	&	1.7485	&	0.4358	&	1.1367	&	1.1650	&	0.1541	&	0.0068	&	0.0251	&	0.9682	\\
14	&	1.7359	&	1.8138	&	0.4511	&	1.1554	&	1.1790	&	0.1564	&	0.0049	&	0.0169	&	0.9782	\\
15	&	1.7975	&	1.8782	&	0.4668	&	1.1733	&	1.1936	&	0.1591	&	0.0035	&	0.0120	&	0.9846	\\
16	&	1.8593	&	1.9419	&	0.4800	&	1.1912	&	1.2086	&	0.1612	&	0.0024	&	0.0082	&	0.9893	\\
17	&	1.9200	&	2.0048	&	0.4938	&	1.2081	&	1.2233	&	0.1629	&	0.0017	&	0.0059	&	0.9924	\\
18	&	1.9802	&	2.0666	&	0.5058	&	1.2252	&	1.2385	&	0.1646	&	0.0012	&	0.0041	&	0.9947	\\
19	&	2.0391	&	2.1273	&	0.5182	&	1.2422	&	1.2536	&	0.1661	&	0.0009	&	0.0029	&	0.9963	\\
20	&	2.0983	&	2.1889	&	0.5313	&	1.2584	&	1.2692	&	0.1673	&	0.0005	&	0.0020	&	0.9974	\\
50	&	3.5612	&	3.6942	&	0.7950	&	1.7125	&	1.7207	&	0.1924	&	-	&	-	&	1	\\
100	&	5.3353	&	5.5101	&	1.0603	&	2.3197	&	2.3295	&	0.2283	&	-	&	-	&	1	
\end{tabular}
}
\caption{Summary statistics for computed $\rho^{\gepp}(A)$ and $\rho^{\gecp}(A)$ using $10^6$ samples  $A \sim \Ginibre(n,n)$}
\label{t:Gauss}
\end{table}

\begin{table}[!htbp]
\centering
{
\begin{tabular}{r|ccccccccc}
\\
$n$& 2&	3&	4&	5&	6&	7&	8&	9&	10	\\ \hline
$\O(n)$&0&	0.2988&	0.5852&	0.8285&	1.0879&	1.2194&	1.416&	1.8546&	1.8546\\
$\GL_n(\mathbb R)$&0&	1&	1.2277&	2.5609&	2.5609&	2.5609&	2.5609&	2.5609&	2.5609\\ \vspace{-.5pc}\\
$n$&11&	12&	13&	14&	15&	16&	17&	18&	19\\\hline
$\O(n)$&	1.8546&	1.9821&	1.9821&	1.9821&	2.0948&	2.3315&	2.3315&	2.4204&	2.8118\\
$\GL_n(\mathbb R)$&	2.5609&3.1294&	3.1294&	3.1294&	3.1294&	3.9719&	3.9719&	3.9719&	3.9719\\ \vspace{-.5pc}\\
$n$&	20&	50&	100\\\hline
$\O(n)$&	3.2711&	5.7837&	7.5449\\
$\GL_n(\mathbb R)$&	3.9719&	7.3208&	11.733
\end{tabular}
}
\caption{Computed lower bounds for $c_n$ for $\O(n)$ and $\GL_n(\mathbb R)$.}
\label{t:c_n}
\end{table}

\begin{table}[!htbp]
\centering
{
\begin{tabular}{r|cccc|cccc}
    &\multicolumn{4}{c}{$X_n = \rho^{\gecp}(Q_n) - \rho^{\gepp}(Q_n+ \overline G)$} &\multicolumn{4}{c}{$Y_n = \rho^{\gecp}(W_n) - \rho^{\gepp}(W_n+ \overline G)$}\\
     $n$ & Median & $\overline x$ & $s$ & \hspace{-.5pc}$\approx \P(X_n = 0)$ & Median & $\overline x$ & $s$ & \hspace{-.5pc}$ \approx \P(Y_n = 0)$\\ \hline 
2	&	2.84e-3	&	0.0032	&	0.0019	&	0.9961	&	2.01e-3	&	0.0023	&	0.0013	&	1.0000	\\
3	&	-9.31e-5	&	-0.4199	&	0.5958	&	0.6669	&	9.89e-4	&	-0.6637	&	0.9404	&	0.6665	\\
4	&	1.00e-5	&	-0.6419	&	1.3045	&	0.7498	&	9.84e-4	&	-0.9030	&	1.8486	&	0.7505	\\
5	&	5.52e-5	&	-0.9290	&	2.4497	&	0.7995	&	9.17e-4	&	-0.9396	&	2.6996	&	0.8001	\\
6	&	7.63e-5	&	-1.4046	&	4.4950	&	0.8330	&	8.49e-4	&	-0.9144	&	3.6066	&	0.8333	\\
7	&	9.50e-5	&	-2.2121	&	8.1828	&	0.8573	&	7.89e-4	&	-0.8642	&	4.6407	&	0.8573	\\
8	&	1.25e-4	&	-3.5869	&	14.747	&	0.8748	&	7.37e-4	&	-0.8128	&	5.9164	&	0.8755	\\
9	&	1.12e-4	&	-5.8481	&	26.175	&	0.8887	&	6.92e-4	&	-0.7732	&	7.6716	&	0.8889	\\
10	&	1.13e-4	&	-9.6106	&	45.923	&	0.9001	&	6.53e-4	&	-0.7222	&	9.5064	&	0.8999	\\
11	&	1.14e-4	&	-15.614	&	78.089	&	0.9091	&	6.19e-4	&	-0.6990	&	12.557	&	0.9092	\\
12	&	1.14e-4	&	-23.934	&	124.76	&	0.9165	&	5.90e-4	&	-0.6236	&	14.360	&	0.9169	\\
13	&	1.12e-4	&	-34.364	&	185.42	&	0.9231	&	5.65e-4	&	-0.5931	&	16.586	&	0.9230	\\
14	&	1.88e-4	&	-46.564	&	255.24	&	0.9281	&	5.41e-4	&	-0.5514	&	19.812	&	0.9284	\\
15	&	1.29e-4	&	-57.685	&	322.47	&	0.9333	&	5.22e-4	&	-0.4736	&	18.366	&	0.9336	\\
16	&	1.08e-4	&	-67.262	&	382.40	&	0.9379	&	5.03e-4	&	-0.4586	&	23.037	&	0.9375	\\
17	&	1.08e-4	&	-76.783	&	440.15	&	0.9411	&	4.87e-4	&	-0.4781	&	24.960	&	0.9413	\\
18	&	1.05e-4	&	-84.588	&	486.28	&	0.9442	&	4.72e-4	&	-0.4483	&	25.173	&	0.9442	\\
19	&	1.05e-4	&	-89.732	&	521.30	&	0.9475	&	4.59e-4	&	-0.3602	&	15.770	&	0.9473	\\
20	&	1.04e-4	&	-94.261	&	549.62	&	0.9502	&	4.46e-4	&	-0.3140	&	13.059	&	0.9507	\\
50	&	1.82e-4	&	-99.268	&	783.61	&	0.9800	&	2.82e-4	&	-0.1059	&	8.9006	&	0.9798	\\
100	&	1.62e-4	&	-74.840	&	817.42	&	0.9900	&	2.06e-4	&	-0.0599	&	13.495	&	0.9900	
\end{tabular}
}
\caption{Summary statistics (median, mean, standard deviation) for $10^6$ samples of $\vep$-neighbors for the difference of initial GECP growth factor to neighbor GEPP growth factor for each of  $Q_n$ and $W_n$, with $\vep = 10^{-3}$, $n=2:20,50,100$, using $\P(|X_n| < \operatorname{tol}) \approx \P(X_n = 0)$ for $\operatorname{tol}=0.01$ and similarly $\P(|Y_n| < \operatorname{tol}) \approx \P(Y_n = 0)$.}
\label{t:Q_n pp-cp}
\end{table}

\section{Maximal orthogonal growth using \texorpdfstring{{\boldmath$Q_n$}}{n}}
\label{sec: max orth prop}

We will show the maximal growth encountered for orthogonal matrices at any intermediate GEPP step is attained by any orthogonal matrix of the form $Q = D Q_n \tilde D$ for sign diagonal matrices $D,\tilde D$, where $Q_n$ is the orthogonal $QR$ factor of $L_n = \V I - \sum_{i>j} \V e_i \V e_j^T$. Again, this supports the conjecture that $g^{\gepp}(\O(n)) = \rho^{\gepp}(Q_n)$, but is insufficient on its own since stronger orthogonal optimization techniques are needed to also account for the denominator in $\rho(A)$. Let $\mathcal Q_n = \{D Q_n \hat D: D,\hat D \in \diag(\{\pm 1\}^n)\}.$

\begin{proposition}\label{prop: max k}
    If $Q \in \O(n)$, then using GEPP,
    \begin{align}
        \max_{i,j,k} |Q^{(k)}_{ij}| \le \sqrt{2 \alpha_{n-1}}
    \end{align}
    where $\alpha_j = 1+\frac23 (4^{j-1} - 1)$ if $j \ge 1$, with equality iff $Q \in \mathcal Q_n$.
\end{proposition}

\begin{proof}
    Let $Q \in \O(n)$. For $L_{ij} = Q^{(j)}_{ij}/Q^{(j)}_{jj}$, then
    \begin{align}
        Q^{(k)}_{ij} &= Q_{ij} + \sum_{\ell = 1}^{\min(i,k) - 1} \left(\sum_{\mathcal J \in \mathcal T(\min(i,k),\ell)}\prod_{e(r,s) \in \mathcal J} (-L_{rs}) \right) Q_{\ell j},
    \end{align}
    where $\mathcal T(n,m)$ is the collection of decreasing paths in $\mathcal K_n$, the complete graph on $n$ vertices, connecting vertices $n > m$ (where $e(i,j)$ designates an edge between vertices $i$ and $j$). Note $|\mathcal T(n,1)| = 2 |\mathcal T(n-1,1)|$ since half of the decreasing paths connecting vertices $n$ to 1  skip vertex $n-1$. Since $|\mathcal T(n,m)| = |\mathcal T(n-m+1,1)|$ and $|\mathcal T(2,1)| = 1$, then 
    \begin{align}
        |\mathcal T(n,m)| = 2^{n-m - 1}.
    \end{align}
    Using the triangle inequality along with the fact $|L_{ij}| \le 1$ for all $i,j$, we have
    \begin{align}
        |Q^{(k)}_{ij}| \le |Q_{ij}| + \sum_{\ell = 1}^{\min(i,k) - 1} 2^{\min(i,k) - \ell - 1} |Q_{\ell j}| = f_{\min(i,k),i}(|Q\V e_j|)
    \end{align}
    where 
    \begin{align}
        f_{r,s}(\V x) &= f_r([x_1,\ldots,x_{r-1},x_s]^T) \qquad \mbox{for} \qquad r \le s, \qquad \mbox{using}\\
        f_r(\V y) &= y_r + \sum_{\ell = 1}^{r-1} 2^{r - \ell - 1} y_\ell
    \end{align}
    for $\V x \in \mathbb R^n$ and $\V y \in \mathbb R^r$. Next, note for $r < s$ and $\V x = |\V x|$, then
    \begin{align}\label{eq: frs max}
        f_n(\V x) \ge f_{s,s}(\V x) = x_s + \sum_{\ell = r+1}^{s-1} 2^{s - \ell - 1} x_\ell + 2^{s - r - 1}\left( x_r + 2  \sum_{\ell = 1}^{r-1} 2^{r-\ell - 1} x_\ell\right) > f_{r,s}(\V x).
    \end{align}
    It follows
    \begin{align}\label{eq: max fik fn}
        \max_{\footnotesize \begin{array}{c} 1 \le i,k \le n\\ \V x \in \mathbb S^{n-1}\end{array}} {f_{\min(i,k),i}(\V x)} = \max_{\V x \in \mathbb S^{n-1}} {f_n(\V x)}.
    \end{align}
    Maximizing the objective function $f_n(\V x)$ given the constraint $\V x \in \mathbb S^{n-1}$ (say, using a Lagrange multiplier method with $g(\V x) = \|\V x\|_2^2 = 1$\footnote{Taking this approach but instead using the constraint $g(\V x) = \|\V x\|_\infty \le 1$ recovers Wilkinson's original upper bound of $2^{n-1}$.}) yields the maximum $f(\hat {\V x}) = \sqrt{2 \alpha_{n-1}}$ with unique maximizer $\hat{\V  x} = \V y/\sqrt{2\alpha_{n-1}} \in \mathbb R^n$ with $y_j = 2^{n-1-j}$ for $j < n$ and $y_n = 1$. This is precisely the last column of the orthogonal matrix $Q_n$  from \Cref{prop:Q}. 
    
    By \Cref{thm: any max Q}, if $Q \in \mathcal Q_n$, then  
    \begin{align}
        \max_{i,j}|Q^{(k)}_{ij}| \le \max_{i,j}|Q^{(n)}_{ij}| = f_n(\hat {\V x}) = \sqrt{2\alpha_{n-1}}.
    \end{align} 
    Conversely, suppose $\tilde Q \in \O(n)$ attains this maximum. Let $m,i,j$ be such that $|\tilde Q^{(m)}_{ij}| = \max_{i,j,k} |\tilde Q^{(k)}_{ij}|$. Let $D$ be a sign matrix such that $D\tilde Q\V e_j = |Q\V e_j|$, and let $\V x = D\tilde Q\V e_j$ (so that $\|\V x\|_2 = 1)$. Then 
    \begin{align}
        f_n(\hat{\V x}) = |\tilde Q^{(m)}_{ij}| \le f_{\min(m,i),i}(\V x) \le f_n(\hat {\V x}).
    \end{align}
    It follows  $i = j = m = n$ and $\V x = \hat{\V x}$ by \eqref{eq: frs max} and \eqref{eq: max fik fn} and the uniqueness of the maximizer in the above computation. Let $\tilde D$ be such that $Q = D\tilde Q\tilde D$ has positive diagonal and last column, where now $\V x = Q\V e_n$. So
    \begin{align}
        |\tilde Q^{(n)}_{nn}| &= Q^{(n)}_{nn} = x_n  + \sum_{\ell = 1}^{n-1} \left(\sum_{\mathcal J \in \mathcal T(n,\ell)} \prod_{e(r,s) \in \mathcal J} (-L_{rs})\right) x_\ell =  f_n(\hat{\V x})
    \end{align}
    so that 
    \begin{align}
         2^{n-1 - \ell} &= |\mathcal T(n,\ell)| = \left|\sum_{\mathcal J \in \mathcal T(n,\ell)}\prod_{e(r,s) \in \mathcal J} (-L_{rs}) \right|   \quad \mbox{for all $\ell < n$,} \quad \mbox{and so}\\
        1&= \prod_{e(r,s) \in \mathcal J} (-L_{rs})  \quad \mbox{for all $\mathcal J \in \mathcal T(n,\ell)$,} \quad \mbox{and hence}\\
        L_{rs} &= -1 \quad \mbox{ for all $r > s$,}
    \end{align}
    i.e., $Q$ has unit lower triangular GENP factor of the form exactly $L_n$. By \Cref{l:qr lu}, then $Q$ is PP and must also be the orthogonal QR factor for $L_n$, so that $Q = Q_n$ from \Cref{prop:Q}. It follows $\tilde Q \in \mathcal Q_n$.

\section{Larger GECP than GEPP growth}
\label{sec: CP > PP}

This is a complementary analysis to that found in \Cref{subsec: extreme}, which focused on models with much larger GEPP growth than GECP.  First, we can define $C \in \O(3)$ by 
\begin{align}\label{eq: max cp 3}
    {C = \begin{bmatrix}
        2/3 & 2/3 & 1/3 \\ -2/3 & 1/3 & 2/3\\ 1/3 & -2/3 & 2/3
    \end{bmatrix}}. 
\end{align}
This is a (scaled) version of the standard matrix that satisfies $\rho^{\gecp}(C) = g_3^{\operatorname{GECP}} = 2.25$.

\begin{figure}[htbp]
  \centering
  \subfloat[$Q(\boldsymbol \t)B_3$]{%
        \includegraphics[width=0.45\textwidth]{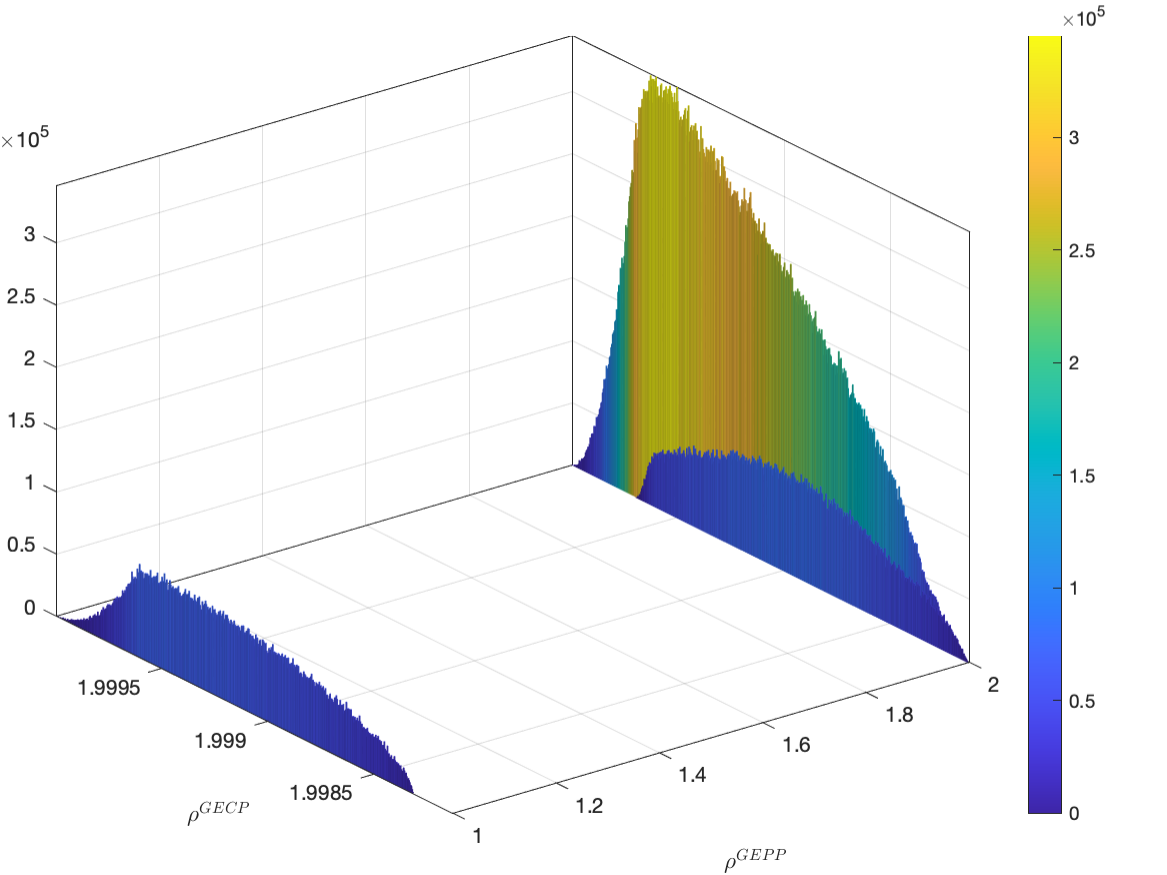}%
        }%
        \label{fig:B1 nbdh1}
    \hspace{-.5pc}%
    \subfloat[$B_3 + (\vep/\sqrt 3)G$]{%
        \includegraphics[width=0.45\textwidth]{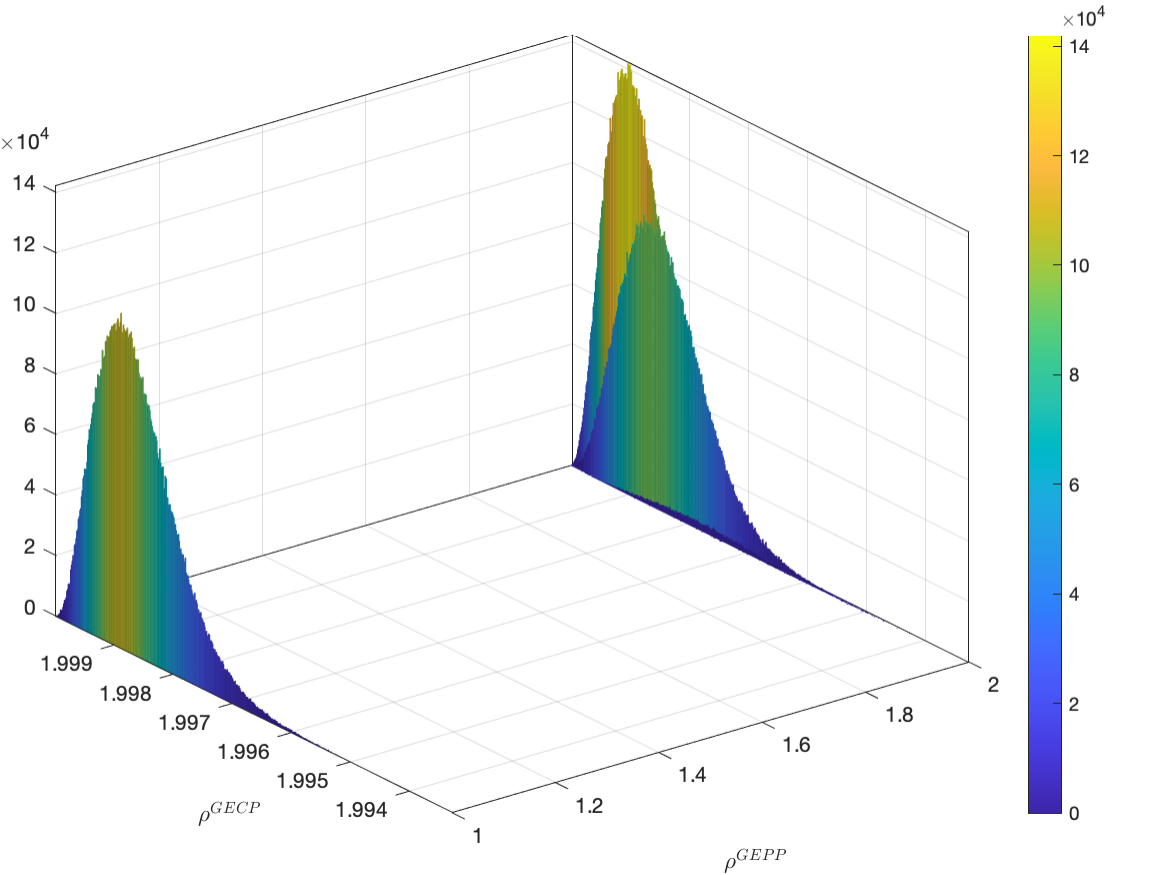}%
        }%
  \caption{Normalized histogram of each pair of $\rho^{\operatorname{GEPP}}$ and $\rho^{\operatorname{GECP}}$ using $10^6$ samples of (a) $Q(\boldsymbol \t)B_3$ for $\boldsymbol \t \sim \Uniform(B_r(\V 0)) \subset \mathbb R^{3}, r = \vep/\sqrt{6}$, and (b) $B_3 + (\vep/\sqrt 3) G$ for $G$ for $G \sim \Ginibre(3,3)$, each using $\vep = 10^{-3}$.}
  \label{fig:B3 nbdh}
\end{figure}

We will first consider  a matrix that has larger GECP growth than GEPP growth. For $n = 3$, the matrix $B_3$ from \eqref{eq: max cp 3} satisfies $2 = \rho^{\gecp}(B_3) > \rho^{\gepp}(B_3) = 1$. We are interested in exploring the behavior for both GEPP and GECP growth on small neighborhoods near $B_3$ (although both remain small since $n$ is so small). \Cref{fig:B3 nbdh} shows a normalized histogram using a $512\times 512$ grid of computed pairs of GEPP and GECP growth factors in $\vep$-neighborhoods near $B_3$ using both orthogonal and additive perturbation methods (see \Cref{subsec: perturb}) for $\vep = 10^{-3}$.

Stability holds for local GECP growth factors, where all computed GECP growth factors for neighbors of $B_3$ stayed within $.005$ of $\rho^{\gecp}(B_3) = 2$. However, the computed GEPP growth factors exhibit highly discontinuous behavior. Each neighbor had GEPP growth that either concentrated only near 1 or 2, and nothing in between.\footnote{All $10^6$ sampled neighbors having GEPP growth factors within 0.002 of either 1 or 2.} This establishes (empirically) both $(2,1)$ and $(2,2)$ comprise limit points for the mapped pairs of growth factors. \Cref{fig:B3 top nbdh} shows the plot of each of the GEPP-GECP growth factor pairs for the neighbors that concentrated near the point $(2,2)$, which comprised the right boundary of \Cref{fig:B3 nbdh}. These nearest neighbors have growth factor pairs that  effectively form an arrow that points toward this limit point $(2,2)$.

\begin{figure}[htbp]
  \centering
  \subfloat[Front view]{%
        \includegraphics[width=0.45\textwidth]{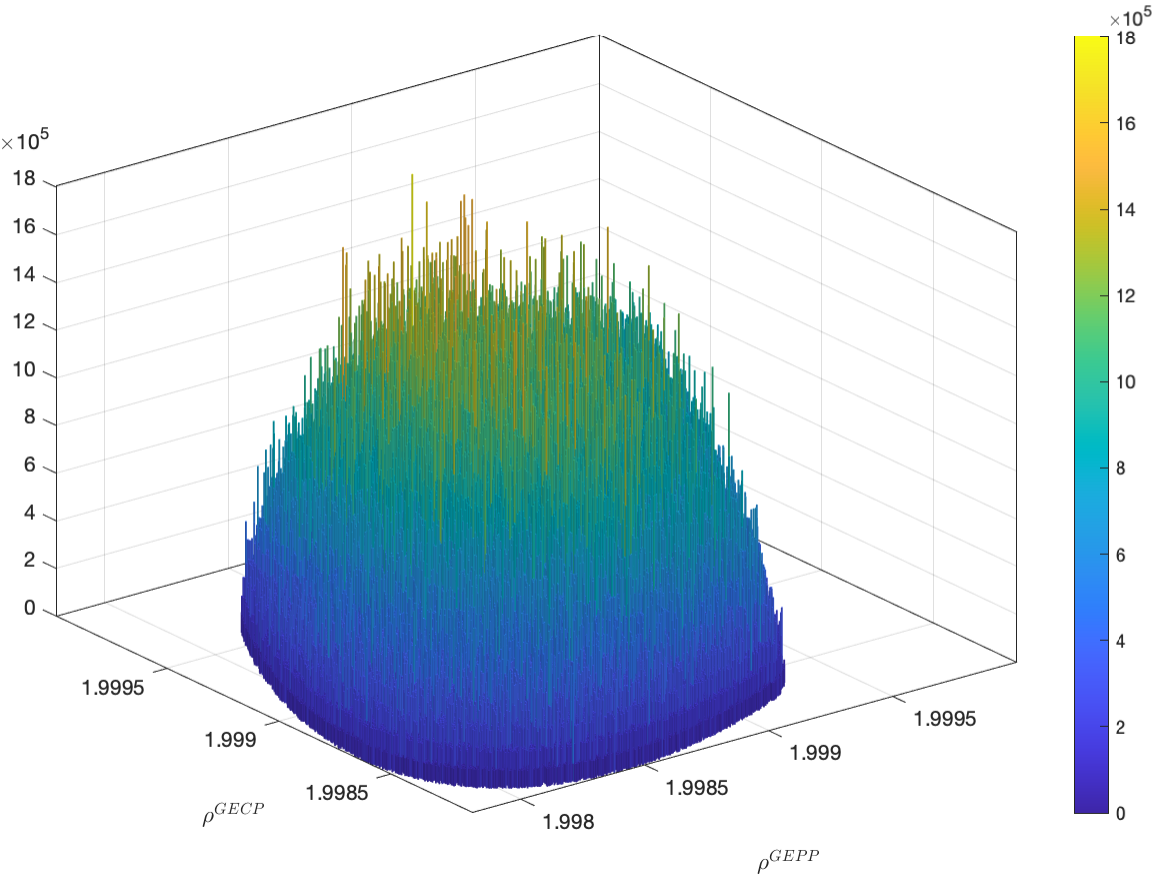}%
        }%
    \subfloat[Top view]{%
        \includegraphics[width=0.45\textwidth]{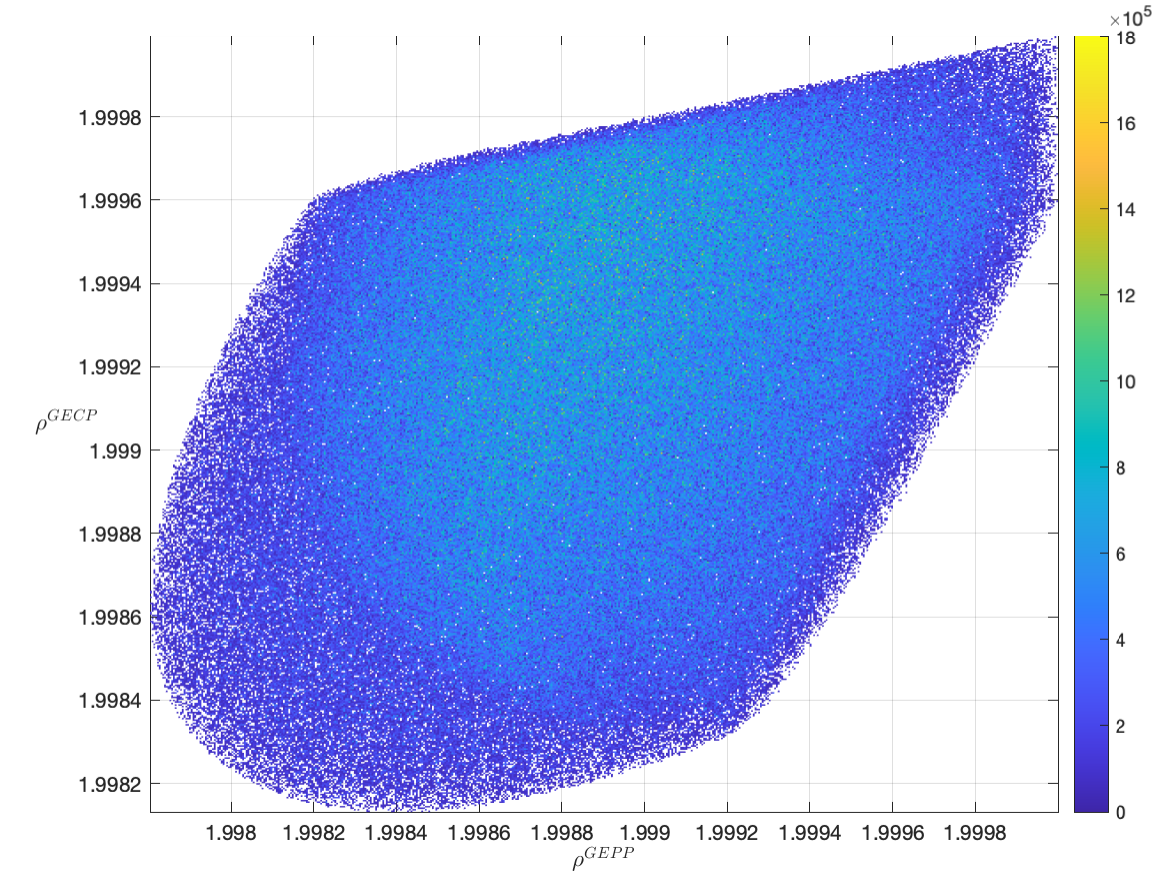}%
        }%
  \caption{Local neighborhood near $(\rho^{\gepp},\rho^{\gecp}) = (2,2)$ using orthogonal perturbations $Q(\boldsymbol \t) B_3$ from \Cref{fig:B3 nbdh}(a).}
  \label{fig:B3 top nbdh}
\end{figure}

What is surprising for this model is that there is higher concentration among these neighbors mapping near $(2,2)$ than the starting point of $(2,1)$. In particular, this shows that while the GECP growth factor remains stable in small neighborhoods, this extreme model has neighbors whose GEPP growth factor moves closer to the \textit{larger} initial GECP factor. For the $10^6$ samples, neighbors were twice as likely to be near $(2,2)$ than $(2,1)$, with 66.68\% of the $10^6$ samples having computed GEPP growth factors within 0.01 of the initial GECP limit point compared to 33.31\% concentrating within 0.01 of the initial GEPP limit point.

\begin{figure}[htbp]
  \centering
  \subfloat[Front view]{%
        \includegraphics[width=0.45\textwidth]{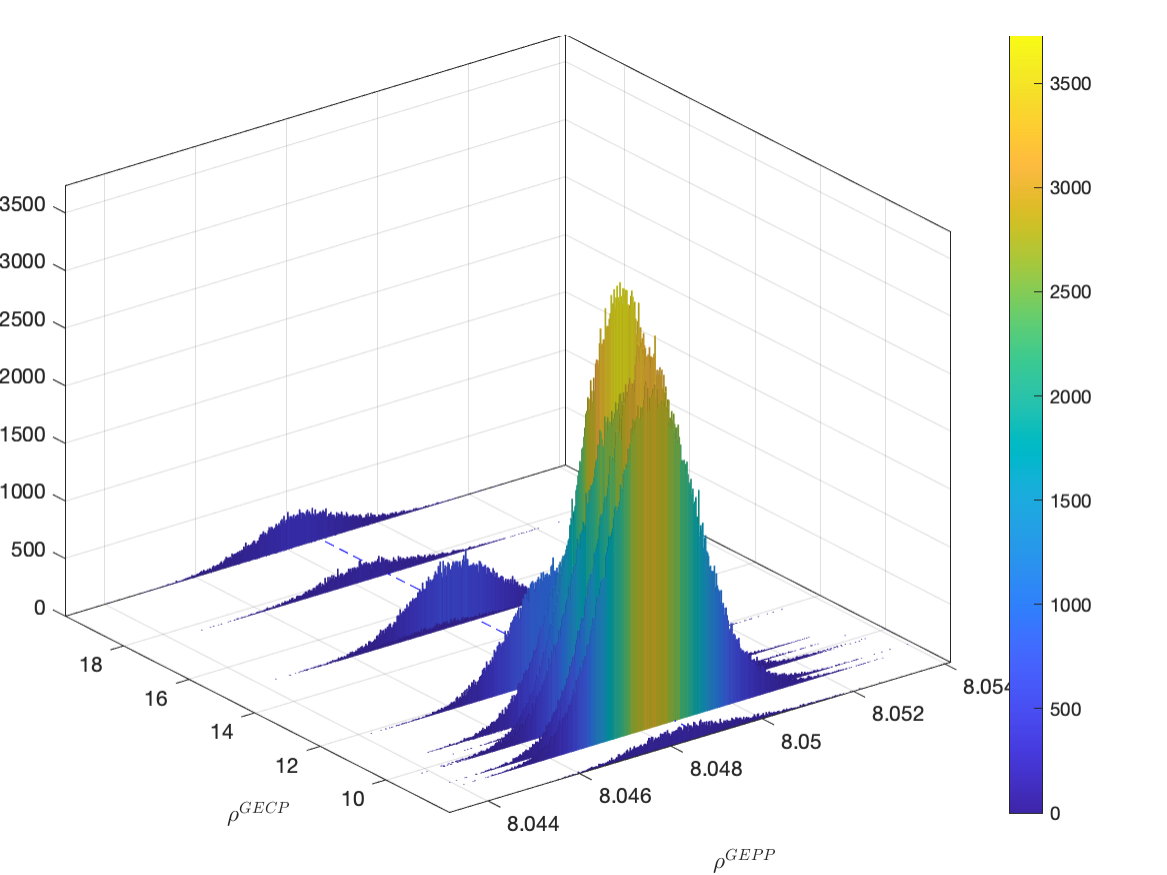}%
        }%
    \subfloat[Top view]{%
        \includegraphics[width=0.45\textwidth]{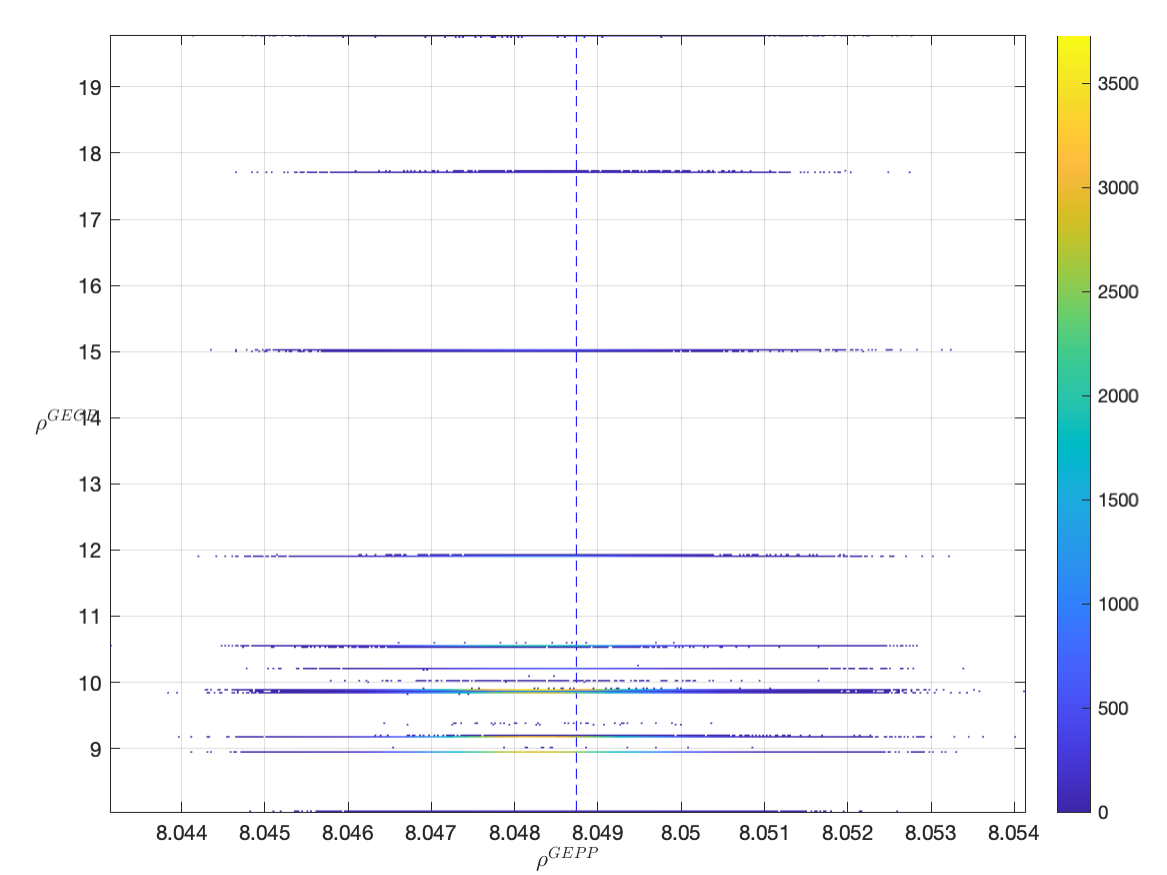}%
        }%
  \caption{Normalized histogram of each pair of $\rho^{\operatorname{GEPP}}$ and $\rho^{\operatorname{GECP}}$ using $10^6$ samples of $B_{100} + (\vep/10) G$ for $G$ for $G \sim \Ginibre(100,100)$ with $\vep = 10^{-3}$, where $B_{100}$ is the computed worst-case matrix that maximizes the GECP-GEPP growth for $n = 100$, with initial $(\rho^{\gepp},\rho^{\gecp}) = (8.0487   ,19.7819)$. A plot of the line $\rho^{\gepp}  = 8.0487$ is included for comparison.}
  \label{fig:B100}
\end{figure}

This behavior does not hold for the remaining computed matrices used for the best lower bounds for $c_n$ from \Cref{t:c_n}. For $n \ge 4$, each of the other computed maximal GECP-GEPP difference models maintained stable GEPP growth factors, with all $10^6$ samples for each $n$ having computed GEPP growth factors stay within 0.02 of the starting GEPP growth factor; so only one limit point held for the GEPP growth factors. However, for these models then the larger initial GECP growth factor showed less stability: as $n$ increased, then the computed GECP growth factors for neighbors moved progressively away from the  larger initial GECP growth factor and closer to the smaller initial GEPP growth factor. \Cref{fig:B100} shows a plot of the computed GEPP and GECP growth factors for $\vep$-neighbors of $B_{100}$, which has the computed maximal GECP-GEPP growth difference for $100\times 100$ matrices of $11.7332$, with $\rho^{\gepp}(B_{100}) = 8.0487$ and $\rho^{\gecp}(B_{100}) = 19.7819$. This model exhibits stable GEPP growth factors but highly discontinuous GECP growth factors that are much more concentrated near 8.0487   than 19.7819. $B_{100}$ also exhibits the first positive proportion  (0.07\%) of the $10^6$ sampled neighboring GECP growth factors appearing within 0.01 of the initial GEPP growth factor among all cases $n = 2:20,50$. 
\end{proof}

\bibliographystyle{siamplain}
\bibliography{references}

\end{document}